\documentclass[12pt]{amsart}
\usepackage{amssymb,times}
\usepackage{amsmath}
\numberwithin{equation}{section}
\def\ca{{\mathcal A}}

\def\ce{{\mathcal E}}
\def\cf{{\mathcal F}}

\def\ch{{\mathcal H}}

\def\car{{\mathcal R}}
\def\cs{{\mathcal S}}

\def\cz{{\mathcal Z}}

\def\ga{{\mathfrak A}}
\def\gb{{\mathfrak B}}
\def\gc{{\mathfrak C}}
\def\gd{{\mathfrak D}}
\def\ge{{\mathfrak E}}

\def\gam{{\mathfrak M}}
\def\gn{{\mathfrak N}}

\def\gar{{\mathfrak R}}

\def\bc{{\mathbb C}}

\def\bm{{\mathbb M}}
\def\bn{{\mathbb N}}

\def\br{{\mathbb R}}
\def\bt{{\mathbb T}}

\def\bz{{\mathbb Z}}

\def\a{\alpha}
\def\b{\beta}
\def\g{\gamma}  \def\G{\Gamma}
\def\d{\delta}  
\def\ve{\varepsilon}

\def\l{\lambda} \def\L{\Lambda}

\def\m{\mu}

\def\n{\nu}
\def\r{\rho}
\def\s{\sigma} 
\def\t{\tau}
\def\f{\varphi} \def\F{\Phi}
  \def\Th{\Theta}
\def\om{\omega} \def\Om{\Omega}

\newtheorem{thm}{Theorem}[section]
\newtheorem{lem}[thm]{Lemma}
\newtheorem{cor}[thm]{Corollary}
\newtheorem{prop}[thm]{Proposition}
\newtheorem{defin}[thm]{Definition}

\def\spec{\mathop{\rm spec}}
\def\id{\mathop{\rm id}}
\def\tr{\mathop{\rm Tr}}
\def\sp{\text{spec}}
\def\di{\mathop{\rm d}\!}

\newcommand{\ty}[1]{\mathop{\rm {#1}}}
\def\idd{{1}\!\!{\rm I}}
\newcommand{\be}{\begin{equation}}
\newcommand{\ee}{\end{equation}}
\newcommand{\bes}{\begin{equation*}}
\newcommand{\ees}{\end{equation*}}
\newcommand{\lb}{\label}
\newcommand{\ds}{\displaystyle}

\newcommand{\nn}{\nonumber}

\begin{document}
\title[Fermi Markov states]
{Fermi Markov states}
\author{Francesco Fidaleo}
\address{Francesco Fidaleo\\
Dipartimento di Matematica\\
II Universit\`{a} di Roma ``Tor Vergata''\\
Via della Ricerca Scientifica, 00133 Roma, Italy}
\email{{\tt fidaleo@mat.uniroma2.it}}

\begin{abstract}
We investigate the structure of the
Markov states on general Fermion algebras. The situation treated in the present paper
covers, beyond the $d$--Markov states on the CAR algebra on $\bz$ (i.e. when there are $d$--annihilators and creators on each site), also the non homogeneous case (i.e. when the numbers of generators depends on the localization). The present analysis provides the first necessary step for the study of the general properties, and the construction of nontrivial examples of Fermi Markov states on 
$\bz^{\n}$, that is the Fermi Markov fields. Natural connections with the KMS boundary condition and entropy of Fermi Markov states are studied 
in detail. Apart from a class of Markov states quite similar to those arising in the tensor product algebras 
(called "strongly even" in the sequel), other interesting examples of Fermi Markov states 
naturally appear. Contrarily to the strongly even examples, the latter are highly entangled and it is expected that they describe interactions which are not  "commuting nearest neighbor". Therefore, the non strongly even Markov states, in addition to the natural applications to quantum statistical mechanics, might be of interest for the information theory as well.
\end{abstract}

\maketitle

\centerline{DEDICATED to L. ACCARDI in OCCASION of his 60th BIRTHDAY}

\section{introduction}

The quantum analogues of Markov processes were first 
constructed in \cite{A}, where the notion of quantum Markov chain on 
infinite tensor product algebras was introduced. 

Nowadays, quantum Markov chains have become a standard computational tool 
in solid state physics,
and several natural applications have emerged in quantum statistical mechanics
and quantum information theory. 
On the other hand, the introduction in \cite{Ar0, AM2, AM3} of the notion of  ``product state'' 
on CAR algebras motivated the analogue construction  in \cite{AFM}, of
quantum Markov chains on these algebras as local perturbation of product states.
The Fermi extension of product and Markov states is nontrivial because, 
even if the Fermion algebra is isomorphic to an infinite tensor products of matrix algebras, this embedding does not 
preserve the natural localization which plays an essential role in the very definition
of these states.
Product states describe non interacting (free or independent) systems. 
Markov states describe nearest neighbor interactions. In both cases the 
notion of "localization" plays a crucial role. 
Typically
a discrete system is identified to a point in a graph. If this graph is not isomorphic to
an interval in $\bz$ ($1$--dimensional case), one speaks of a random field.
The crucial role of the localization is at the root of the difficulties to 
construct nontrivial examples of Markov fields. As the ``interacting degrees of freedom''  localized in a finite volume, increase with the volume, the first step to achieve 
this goal is to investigate the nonhomogeneous $1$--dimensional case, one of the goals of the present paper. Our second goal has to do with the most important difference
between tensor and Fermi Markov chains, emerged from the analysis of \cite{AFM},
which in its turn is related to 
the difference between {\it quantum Markov chains} and {\it quantum Markov states}.
The origin of this difference lies in the fact that, in the classical case,
the simple structure of Markov states is equivalent to a single intrinsic condition: 
the Markov property. In the quantum case, while the Markov property can be formulated in terms of a  localization property of the modular group 
of the state (see \cite{AF2, AFr}), there is a class of states which have a 
Markov like local structure but do not necessarily enjoy the Markov property.
This states are called quantum Markov chains (see e.g. \cite{AFM}, Definition 2.2).
This phenomenon is related to the fact that the 
natural probabilistic extension of the notion of conditional expectation to
the quantum case in general is not a projection (cf. \cite{AC}).
On the other hand, the results of \cite{FNW} show that 
some of the most interesting physical applications involve 
precisely those Markov chains which are not Markov states. Such Markov chains can be explicitly constructed, but at the moment no intrinsic operator theoretic 
characterization is known.
This distinction also appears in the Fermi case. However, the new phenomenon consists in the fact that, while in the tensor case
all Markov states are convex combinations of states which are product states with
respect to a new localization canonically associated to the original one
(two block factors), this is not true in the Fermi case. 
The Fermi analogue of the convex combination of two block factors still appear. The last are called in the sequel {\it strongly even} Markov states.
In addition, there is a completely new class of Fermi Markov states.
This new class of non strongly even Markov states that appears in the classification theorem below
(cf. Theorem \ref{dis}) is likely to play in the Fermi case, the role played by the entangled states
in the tensor product case. 

One can define the
notion of Fermi entanglement in analogy to the tensor product case.
One can expect that the main problem of the entanglement theory, that is to
find constructive and easily applicable criteria to discriminate entangled from
non entangled states, will be in the Fermi case at least as difficult as in the tensor case.
The first step to attack this problem, is to have a full and detailed description of this new class of states. In the paper \cite{AFM} only the simplest case was considered, that is when there is only one creator and annihilator in each site. 
The non homogeneous case, discussed in the present paper, includes as a particular case
the translation invariant cases described in \cite{AFM} and its natural translation invariant generalization when there are $d$ creators and annihilators localized in each site. This leads to a much larger class of non strongly even Markov states, which can be completely described. 

In the present paper, the program outlined below is carried out in the following steps.  
Section 2 contains the key result on the structure of the even transition 
expectations associated to Fermi Markov states (Proposition \ref{part}).
It is then possible to provide the full classification of all even Markov quasi conditional expectations.
Section 3 is devoted to the study of the most general situation, including also the non
homogeneous Fermi Markov states. Even if the structure of the Markov states considered here is more complex than the one in \cite{AF1, AFM}, we are still able to provide their decomposition as  direct integrals of
minimal ones (Theorem \ref{dis}). Furthermore, the minimal Markov states are the building--blocks for the construction of all the Fermi Markov states (Theorem \ref{disbis}).  
Section 4 deals with some general properties of the Markov states, 
such as the connection with the KMS boundary condition, and the entropy. 
Section 5 is devoted to the study of the detailed structure of the strongly even Markov states. They can be viewed as the Fermi analogue of the Ising type interactions. We show 
that a strongly even Markov state $\f$ on the Fermion algebra arises by a lifting of a classical Markov process on the spectrum of a maximal Abelian subalgebra, with respect to the same localization as $\f$. In addition, we establish the equality between the Connes--Narnhofer--Thirring dynamical entropy 
$h_{\f}(\a)$ with respect the shift and the mean entropy $s(\f)$.
Section 6 provides the full list of translation invariant 
Fermi Markov states for low dimensional single site local algebras.
The same method is applicable to higher dimensions. Finally, by using Moriya criterion, we show (cf. Proposition \ref{eeeent})  that the Fermi Markov states which are not strongly even are indeed entangled, that is they provide a wide class of examples of entangled states on the CAR algebra which can be directly constructed and investigated in detail.

\section{preliminaries}

For the convenience of the reader, we collect some preliminary facts needed in the 
sequel.

\subsection{Umegaki conditional expectations}

By a (Umegaki) {\it conditional expectation} 
$E:\ga\to\gb\subset\ga$ we mean a norm one projection of the 
$C^{*}$--algebra $\ga$ onto the $C^{*}$--subalgebra (with the same 
identity $\idd$) $\gb$.  The case of interest for us is 
when $\ga$ 
is a full matrix algebra. Consider a set 
$\{P_{i}\}$ of central orthogonal projections of the range $\gb$ of 
$E$,
summing up to the 
identity. 
We have 
\begin{equation}
\label{mmppr}
E(x)=\sum_{i}E(P_{i}xP_{i})P_{i}\,.
\end{equation} 
Then 
$E$ is uniquely determined by its values on
the reduced algebras
$\ga_{P_{i}}:=P_{i}\ga P_{i}$.
When the above set $\{P_{i}\}$ consists of minimal projections, we get
$\ga_{P_{i}}=N_{i}\otimes\bar N_{i}$, 
and there exist states $\phi_{i}$ on $\bar N_{i}$ such that
\begin{equation}
\label{cond}
E(P_{i}(a\otimes\bar a)P_{i})=\phi_{i}(\bar a)P_{i}(a\otimes\idd)P_{i}\,.
\end{equation}
The reader is referred to \cite{H} for further details.
\subsection{Quasi conditional expectations}

Consider a triplet $\gc\subset\gb\subset\ga$ of  
unital $C^{*}$--algebras. A {\it quasi conditional expectation} w.r.t. the given
triplet, is a  completely positive, identity preserving linear map 
$E:\ga\to\gb$ such that
\bes
E(ca)=cE(a)\,, \quad a\in\ga\,,\, c\in\gc\,.
\ees
Notice that, as 
the quasi conditional expectation $E$ is a real map, we have
$$
E(ac)=E(a)c\,, \quad a\in\ga\,,\, c\in\gc\,.
$$
If $\f$ is a normal faithful state on the $W^{*}$--algebra $\ga$,    
the $\f$--expectation $E^{\f}:\ga\to\gb$ by Accardi and Cecchini 
preserving the 
restriction of $\f$ to the $W^{*}$--subalgebra $\gb$,
provides an 
example of quasi conditional expectation. Namely, it is enough to choose for $\gc$ 
any unital $C^{*}$--subalgebra of $\gb$ contained in the 
$E^{\f}$--fixed point algebra. $E^{\f}$ 
is a conditional expectation if and only if the modular group of $\f$ 
leaves globally stable the subalgebra $\gb$, see \cite{AC}.

\subsection{The CAR algebra}

Denote $[a,b]:=ab-ba$, $\{a,b\}:=ab+ba$
the commutator and anticommutator 
between elements $a$, $b$ of an algebra, respectively.

Let $J$ be a set. The {\it Canonical Anticommutation Relations} (CAR for short) algebra 
over  $J$ is the $C^{*}$--algebra $\ga_J$ with the identity $\idd$
generated by the set $\{a_j, a^+_j\}_{j\in I}$, and the relations
\bes
(a_{j})^{*}=a^{+}_{j}\,,\,\,\{a^{+}_{j},a_{k}\}=\d_{jk}\idd\,,\,\,
\{a_{j},a_{k}\}=\{a^{+}_{j},a^{+}_{k}\}=0\,,\,\,j,k\in J\,.
\ees
When there is no matter of confusion, we denote $\ga_J$ simply as $\ga$.
The parity automorphism $\Th$, of $\ga$ acts on the generators as
$$
\Th(a_{j})=-a_{j}\,,\,\,\Th(a^{+}_{j})=-a^{+}_{j}\,,\quad j\in J\,, 
$$
and induces on $\ga$ the $\bz_{2}$--grading
$\ga=\ga_{+} \oplus \ga_{-}$ where
$$
\ga_{+}:=\{a\in\ga \ | \ \Th(a)=a\}\,,\quad
\ga_{-}:=\{a\in\ga \ | \ \Th(a)=-a\}\,.
$$
Elements  in $\ga_+$ (resp. $\ga_-$) are called 
{\it even} (resp. {\it odd}).

A map $T:\ga^1\to\ga^2$ between the CAR algebras $\ga^1$, $\ga^2$  with $\bz_{2}$--gradings
$\Th_1$, $\Th_2$ is said to be 
{\it even} if it 
is grading--equivariant: 
$$
T\circ\Th_1=\Th_2\circ T\,.
$$
The previous definition applied to states $\f\in\cs(\ga)$ leads to $\f\circ\Th=\f$, that is $\f$ is even if it is
$\Th$--invariant. 

Let the index set $J$ be countable, then the CAR algebra  is isomorphic to the 
$C^{*}$--infinite 
tensor product of $J$--copies of $\bm_{2}(\bc)$:
\begin{equation}
\label{jkw}
\ga_J\sim\overline{\bigotimes_{J}\bm_{2}(\bc)}^{C^*}\,.
\end{equation}
For the convenience of the reader, we report the Jordan--Klein--Wigner transformation establishing the mentioned 
isomorphism. Fix any enumeration $j=1,2,\dots$ of the set $J$. Let
$U_{j}:=a_{j}a_{j}^{+}-a_{j}^{+}a_{j}$, $j=1,2,\dots$\,. Put $V_{0}:=\idd$, 
${\ds V_{j}:=\prod_{n=1}^{j}U_{n}}$, and denote
\begin{align}
\lb{kw}
&e_{11}(j):=a_{j}a_{j}^{+}\,,\quad e_{12}(j):=V_{j-1}a_{j}\,,\nn\\
&e_{21}(j):=V_{j-1}a_{j}^{+}\,,\quad e_{22}(j):=a_{j}^{+}a_{j}\,.
\end{align}
$\{e_{kl}(j)\,|\,k,l=1,2\}_{j\in I}$ 
provides a system of commuting $2\times2$ matrix units realizing the 
mentioned isomorphism. 

Thanks to \eqref{jkw}, $\ga_J$ has a unique tracial state $\t$ (at least when $J$ is countable) as the extension of the unique tracial state on $\ga_I$, $|I|<+\infty$. 
Let $J_1\subset J$ be a finite set and
$\f\in\cs(\ga)$. Then there exists a unique positive element $T$ such that 
$\f\lceil_{\ga_{J_1}}=\t\lceil_{\ga_{J_1}}(\,{\bf\cdot}\,T)$. The element
$T$ is called the {\it adjusted matrix} of $\f\lceil_{\ga_{J_1}}$.\footnote{For the standard applications to quantum statistical mechanics, we also use the density matrix w.r.t. the unnormalized trace, see Section 5.} 
The state $\f\lceil_{\ga_{J_1}}$ is even (faithful) if and only if its adjusted matrix is even (invertible). 
The reader is referred to \cite{T}, Section XIV.1 and \cite{AM2} for further details.

We end the present subsection by recalling the description of product state (cf. \cite{AM2}), and the definition of entanglement (cf. \cite{M}, Section 2). 
Let $J_1,J_2\subset I$ with $J_1\bigcap J_2=\emptyset$. Fix $\f_{1}\in\cs(\ga_{J_1})$, 
$\f_{2}\in\cs(\ga_{J_2})$. If at least one among them is even, then according to Theorem 11.2 of 
\cite{AM2}, the product state extension (called {\it product state} for short) 
$\f\in\cs(\ga_{J_1\bigcup J_2})$ is uniquely defined. We write with an abuse of notation, $\f=\f_1\f_2$. Suppose that $J_1,J_2$ are finite sets. Let 
$T_{1}\in\ga_{J_1}$, 
$T_{2}\in\ga_{J_2}$ be the adjusted densities
relative to $\f_{1}\in\cs(\ga_{J_1})$, $\f_{2}\in\cs(\ga_{J_2})$, respectively. If at least one among $T_{1}$ and $T_{2}$ is even, then $[T_{1},T_{2}]=0$ and $T:=T_{1}T_{2}$ is a well defined 
positive element of $\ga_{J_1\cup J_2}$ which is precisely 
the density matrix of $\f=\f_1\f_2$. $\f\in\cs(\ga_{J_1\cup J_2})$ is even if and only if $\f_1$ and $\f_2$ are both even.

A state $\f\in\cs(\ga_{J_1\cup J_2})$ is called {\it separable} (w.r.t. to the decomposition
$\ga_{J_1\cup J_2}=\overline{\ga_{J_1}\bigvee\ga_{J_2}}$) if it is in 
the closed convex hull of all the product states over $\ga_{J_1\cup J_2}$. Otherwise it is called {\it entangled}.

\subsection{Preliminaries on Fermi Markov states}

Let us start as in \cite{AF1},  
with a totally ordered countable set $I$
containing, possibly a smallest element $j_{-}$ and/or a greatest
element $j_{+}$. If $I$ contains neither $j_{-}$, nor
$j_{+}$, then $I\sim\bz$. If only $j_{+}\in I$, then $I\sim\bz_{-}$, and
if only $j_{-}\in I$, then $I\sim\bz_{+}$. Finally, if both
$j_{-}$ and $j_{+}$ belong to $I$, then $I$ is a finite set and the 
analysis becomes easier. If  
$I$ is order isomorphic to $\bz$, $\bz_{-}$ or $\bz_{+}$, we put
simbolically $j_{-}$ and/or $j_{+}$ equal to $-\infty$ and/or $+\infty$
respectively. In such a way, the objects with indices
$j_{-}$ and $j_{+}$ will be missing in the computations.

Let $\ga_j$ be the CAR algebra generated by $d_j$ creators and annihilators 
$\{a_{j,1},a_{j,1}^{+},
a_{j,2},a_{j,2}^{+},\dots,a_{j,d_{j}},a_{j,d_{j}}^{+}\}$ localized on the site 
$j\in I$. The numbers of  the $2d_j$ generators of $\ga_j$ may depend on $j$. 
We call 
\begin{equation}
\label{feal}
\ga:=\overline{\bigvee_{j\in I}\ga_{j}}^{C^*}
\end{equation}
{\it the Fermion algebra}. Let ${\displaystyle J:=\bigcup_{j\in I}\{1,2,\dots d_j\}}$ be the disjoint union of the sets $\{1,2,\dots d_j\}$, $j\in I$. Then the Fermion algebra $\ga$ given in \eqref{feal} is nothing but the CAR algebra over the set $J$ previously described.

Now we pass to describe the local structure of the Fermion algebra $\ga$. For each $\L\subset I$, the local algebra $\ga_{\L}\subset\ga$ is defined as 
${\displaystyle\ga_{\L}:=\overline{\bigvee_{j\in\L}\ga_j}}$. According to this notation, $\ga_{\{j\}}=\ga_{j}$ and $\ga_I=\ga$. Then $\L\subset I\mapsto\ga_{\L}\subset\ga$ describes the local structure of the Fermion algebra. Particular subsets of $I$ are
$$
[k,n]:=\big\{l\in I\,\big|\,k\leq l\leq n\big\}\,,\quad
n]:=\big\{l\in I\,\big|\,l\leq n\big\}\,.	
$$
We put for $\L\subset I$, $S_{\L}:=S\lceil_{\ga_{\L}}$, $S$ being any map defined  
on $\ga$. The reader is referred to \cite{BR1}, Section 2.6 and \cite{AM2}, Section 4 for further details. 

A state $\f\in\ga$ is said to be {\it locally faithful} if $\f_{\L}$ is faithful whenever $\L\subset I$ is finite. 

If the number local generators $d_j$ depend on $j$ we refer to this situation as the {\it nonhomogeneous case}. Conversely, when $I=\bz$ and $d_j=d$, $j\in\bz$, the shift  
$j\longrightarrow j+1$ acts in a natural way as an automorphism $\a$ of $\ga$. A state $\f\in\cs(\ga)$ is {\it translation invariant} if $\f\circ\a=\f$. If a state is translation invariant, then it is automatically even, 
see e.g. \cite{BR1}, Example 5.2.21.

We pass to the definition of Markov states which parallels Definition 4.1 of \cite{AFM}.
\begin{defin}
\lb{dema}
An even state $\f$ on $\ga$ is called a Markov state
if, for each $n<j_{+}$, there exists an even quasi
conditional expectation
$E_{n}$ w.r.t. the triplet 
$\ga_{n-1]}\subset\ga_{n]}\subset\ga_{n+1]}$ satisfying
\begin{equation}
\lb{mrec}
\f_{n]}\circ E_{n}=\f_{n+1]}\,,
\end{equation}
$$
E_{n}(\ga_{[n,n+1]})\subset\ga_{\{n\}}\,.
$$
\end{defin}

Notice that the local structure $\L\mapsto\ga_{\L}$, $\L$ finite subset of $I$, plays a crucial role in defining the Markov property. In fact, the isomorphism in \eqref{jkw} does not preserve neither the grading nor the natural localization.\footnote{The algebra on the r.h.s. of \eqref{jkw} is naturally equipped with the trivial parity automorphism. Thus, its $\bz_2$--grading is trivial.} Hence, 
it does not intertwine the corresponding Markov states. 

When the numbers $d_j$ of the generators of 
$\ga_j$ depend of the site, we call a Markov state $\f$ (or equally well a Markov measure in the Abelian case) a {\it nonhomogeneous Markov state}. If $d_j=d$ for each $j$ we refer to the 
{\it $d$--Markov property}. Thus, {\it homogeneity} means $d$--Markov property for some $d$. For the applications to quantum statistical mechanics, $d_j$  is nothing but the "range of interaction" on the chain which might depend on the site, and when $d=1$ we are speaking of {\it nearest neighbor interaction}. The reader is referred to 
\cite{AF1, AF2, Pr} and the literature cited therein, for the connection between the Markov property and the statistical mechanics, and for further details.

Let $\f\in\cs(\ga)$ be a locally faithful Markov state. Then the restriction 
$e_{n}:=E_{n}\lceil_{\ga_{[n,n+1]}}$ is a completely positive identity preserving linear map 
$e_{n}:\ga_{[n,n+1]}\to\ga_{\{n\}}\subset\ga_{[n,n+1]}$ leaving invariant the faithful state 
$\f_{[n,n+1]}$. It is a quite standard fact (see e.g. \cite{AC}) that the ergodic average
\begin{equation*}
\ve_{n}:=\lim_k\frac{1}{k}\sum^{k-1}_{h=0}(e_{n})^{h}
\end{equation*}
exists and defines a conditional expectation 
$$
\ve_{n}:\ga_{[n,n+1]}\to\car(\ve_n)\subset\ga_{\{n\}}
$$
projecting onto the fixed point algebra of $e_{n}$, the last coinciding with the range $\car(\ve_n)$ of 
$\ve_n$.
The sequence $\{\ve_{n}\}_{n<j_+}$ of two point conditional expectations is called in the sequel the sequence of {\it transition expectations} associated to the locally faithful Markov state $\f$. They uniquely determine, and are determined by 
the conditional expectations $\ce_{n}:\ga_{n+1]}\to\ga_{n]}$, given 
for $x\in\ga_{n-1]}$, $y\in\ga_{[n,n+1]}$ by
\begin{equation}
\label{cecp}
\ce_{n}(xy)=x\ve_{n}(y)\,.
\end{equation}

In addition, it is quite standard to verify (cf. \cite{AFM}, Proposition 4.2) that we can freely replace the quasi conditional expectation $E_n$ in Definition \ref{dema} with its ergodic average $\ce_n$. For the convenience of the reader we report Proposition 4.3 of \cite{AFM}.
\begin{prop}
\lb{twst10}
Let $f:\ga_{[n,n+1]}\to\car(f)\subset\ga_{\{n\}}$ be 
a even conditional expectation. The formula
\bes
\cf(xy):=xf(y)\,,\quad x\in\ga_{n-1]}\,,\quad y\in\ga_{[n,n+1]}
\ees
uniquely defines a even conditional expectation 
$$
\cf:\ga_{n+1]}\to\ga_{n-1]}\bigvee\car(f)\subset\ga_{n]}\,.
$$
\end{prop}
From now on, we deal without further mention with even (quasi) conditional expectations. In addition, all the Markov states we deal with are even, and locally faithful if it is not otherwise specified.
\begin{lem}
\lb{fcex}
Let $\ce:\ga_{[k,l+1]}\to\car(\ce)\subset\ga_{[k,l]}$ be a  
conditional expectation with $\ga_{[k,l-1]}\subset\car(\ce)$. 
Then $\ce$ is faithful 
provided that $\ce\big\lceil_{\ga_{[l,l+1]}}$ is faithful.
\end{lem}
\begin{proof}
Let $\f_{1}$, $\psi$ be faithful even states on $\ga_{[k,l-1]}$, 
$\car(\ce\lceil_{\ga_{[l,l+1]}})$ respectively. 
Put $\f_{2}:=\psi\circ\ce\lceil_{\ga_{[l,l+1]}}$.
The product state $\f:=\f_{1}\f_{2}$ 
is a faithful state on $\ga_{[k,l+1]}$ left invariant by 
$\ce$. Fix $a\in\ga_{[k,l+1]}$ with 
$\ce(a^{*}a)=0$. Then $\f(a^{*}a)=\f(\ce(a^{*}a))=0$ which implies 
that $a=0$ as $\f$ is faithful. Namely, $\ce$ is
faithful.
\end{proof}
We then pass to study the structure of the even conditional expectations
$$
\ve_n:\ga_{[n,n+1]}\to\car(\ve_n)\subset\ga_{\{n\}}\,.
$$
To shorten the notations, it is enough to consider the case when $n=0$. After putting $\ve:=\ve_0$, let us start with the finite set $\{P_{j}\}$ of the minimal projections of the 
centre $\cz(\car(\ve))$ of $\car(\ve)$.
\begin{lem}
\lb{orb}
The parity automorphism $\Th$ acts on $\cz(\car(\ve))$, and the orbits 
of minimal projections consist of one or two elements.
\end{lem}
\begin{proof}
Let $P_{j}$ be a minimal projection of $\cz(\car(\ve))$. As $\ve$ is 
even and $\Th^{2}=\id$, we have that $\Th(P_{j})$ is a minimal 
projection of $\cz(\car(\ve))$. This means that either $\Th(P_{j})=P_{j}$, or
$\Th(P_{j})$ is orthogonal to $P_{j}$. The latter means that the 
orbit of $\Th(P_{j})$ consists of two elements.
\end{proof}
We showed in \cite{AFM} that there are interesting examples with 
$\Th(P_{j})\neq P_{j}$.
Let $\ve$ be as above. Some useful properties of the pieces $\ve(PxP)P$, $P$ 
being a even projection of the centre of $\car(\ve)$, minimal among 
the invariant ones, are
described below. 
\begin{lem}
\lb{bos}
Let $M=M_+\oplus M_-$ be a $\bz_{2}$--graded full matrix algebra.
If $x\in M_{-}$ commutes with $M_{+}$, then $x=0$.
\end{lem}
\begin{proof}
Let the $\bz_{2}$--grading be implemented by the automorphism $\Th$.
As $\Th\lceil_{M}$ is inner, there exists an even selfadjoint unitary $V\in M$, 
uniquely determined up to a sign, 
implementing $\Th$ on $M$, see \cite{SZ}, Corollary 8.11. 
This means that $M_{+}=A'$, $A$ being the 
Abelian algebra generated by $V$, and the commutant is taken in the 
full matrix algebra $M$. As $x\in(M_{+})'$, $x\in A''\equiv A$. 
As $x$ is odd, we have $VxV=-x$. Collecting together, we obtain 
$x=0$.
\end{proof}
\begin{lem}
\label{bboss}
Let $M=M_+\oplus M_-$ be a $\bz_{2}$--graded full matrix algebra. For every $\Th$--invariant full matrix subalgebra $N\subset M$, there exists a unique $\Th$--invariant full matrix subalgebra 
$\bar N\subset M$ such that $N\bigvee\bar N$, and
\begin{equation}
\label{ckre}
x\bar x+\s(x,\bar x)\bar x x=0\,,\quad x\in N_{\pm}\,, \bar x\in\bar N_{\pm}
\end{equation}
where $\s(x,\bar x)$ is 1 if both $x$, $\bar x$ are odd, and $-1$ in 
the remaining cases. Moreover, we have $N\bigwedge\bar N=\bc\idd$.
\end{lem}
\begin{proof}
Let $\tilde N:=N'\bigwedge M$ which is a $\Th$--invariant full matrix subalgebra of $M$ as well. Fix a (even) unitary $V\in N$ uniquely determined up to a sign, 
implementing $\Th$ on $N$. Define for $x=x_++x_-\in\tilde N$, 
\begin{equation}
\label{bbeta}
\b(x):=x_++Vx_-\,.
\end{equation}
It is easy to see that 
$\b$ defines a $*$--algebra isomorphism between $\tilde N$ and $\b(\tilde N)$. Thus, the full matrix algebra $\bar N:=\b(\tilde N)$ is the algebra we are looking for. 

For the uniqueness, let $\bar R\subset M$ be a $\Th$--invariant full matrix algebra fulfilling the commutation relations in \eqref{ckre} whenever $x\in N_{\pm}$, $\bar x\in\bar R_{\pm}$, such that 
$N\bigvee\bar R=M$. Then it is easy to verify that 
$\tilde R:=\big\{x_++Vx_-\,\big|\,x\in\bar R\big\}$ is a full matrix subalgebra of $\tilde N$. Since 
$M=N\bigvee\tilde N\sim N\otimes\tilde N$, we get that $\tilde R$ must coincide with $\tilde N$ which implies $\bar R=\bar N$.
\end{proof}
We call the algebra 
\begin{equation}
\label{bbar}
\bar N:=\big(N'\bigwedge M\big)_{+}+V\big(N'\bigwedge M\big)_{-}
\end{equation}
obtained in Lemma \ref{bboss},
the {\it Fermion complement} of $N$ in $M$.

\begin{prop}
\lb{part}
Let $\ga:=\bigvee_{j\in I}\ga_j$ be the Fermi $C^*$--algebra with $I=\{0,1\}$.\\
(i) Let $P\in\ga_{\{0\}}$ be a $\Th$--invariant projection. Then there is a one--to--one correspondence between:
\begin{itemize}
\item[(a)] $\ve:P\ga P\to P\ga P$ an even conditional expectation such that $\car(\ve)$ is a full matrix subalgebra of $P\ga_{\{0\}}P$,
\item[(b)] $N\subset P\ga_{\{0\}}P$ a $\Th$--invariant full matrix subalgebra and $\F$ an even state on 
$\bar N\bigvee P\ga_{\{1\}}P$.
\end{itemize}
The correspondence is given for $x\in N$, $y\in\bar N\bigvee P\ga_{\{1\}}P$ by 
\be
\lb{stpar}
\ve(xy)=
\F(y)x
\ee
where $\bar N$ is the Fermion complement of $N$ in $P\ga_{\{0\}}P$ given in \eqref{bbar}. In particular, $\car(\ve)=N$.\\
\noindent
(ii) Let $P_1,P_2\in\ga_{\{0\}}$ such that $\Th(P_1)=P_2$, $P_1P_2=0$. Then there is a one--to--one correspondence between:
\begin{itemize}
\item[(a)] $\ve:(P_1+P_2)\ga(P_1+P_2)\to (P_1+P_2)\ga(P_1+P_2)$ an even conditional expectation such that $\car(\ve)\subset\ga_{\{0\}}$ and $\cz(\car(\ve))=\bc P_1\oplus\bc P_2$,
\item[(b)] $N_1\subset P_1\ga_{\{0\}}P_1$ full matrix algebra and $\F$ a state on 
$M_1:=N_1'\bigwedge P_1\ga P_1$. 
\end{itemize}
The correspondence is given for $x_{i}\in N_{i}$, $y_{i}\in M_{i}$, $i=1,2$, 
\be
\lb{stpar1}
\ve(x_{1}y_{1}+x_{2}y_{2})=
\F(y_{1})x_{1}+\F(\Th(y_{2}))x_{2}
\ee
where $N_2:=\Th(N_1)$, $M_2:=\Th(M_1)$. In particular, $\car(\ve)=N_1\oplus \Th(N_1)$.

In addition, if $z\in\ga_{\{1\}}$ is even, then
\be
\lb{stpar2}
\ve((P_{1}+P_{2})z)=\F(P_{1}zP_{1})(P_{1}+P_{2})\,.
\ee
\end{prop}
\begin{proof}
(i) Let $N:=\car(\ve)$. As $\ve$ is even, $N$ is a $\Th$--invariant full matrix algebra of $P\ga_{\{0\}}P$.  Let 
$\bar N$ be the Fermion complement of $N$ in $P\ga_{\{0\}}P$,
and $y\in\bar N\bigvee P\ga_{\{1\}}P$
an odd element. Then 
$\ve(y)\in N$ is odd too, and by the bimodule property of $\ve$, 
$\big[\ve(y),N_{+}\big]=0$. By Lemma \ref{bos},
$\ve(y)=0$. If $x\in N$, $y\in\bar N\bigvee\ga_{\{1\}}$, we have
\bes
x\ve(y)=x\ve(y_{+})
=\ve(xy_{+})=\ve(y_{+}x)
=\ve(y_{+})x=\ve(y)x\,.
\ees
This means that $\ve(y)\in\cz(N)\equiv\bc P$, that is 
$\ve(xy)=\F(y)x$ for a uniquely determined even state 
$\F$ on $\bar N\bigvee P\ga_{\{1\}}P$. 

Fix now an invariant full matrix subalgebra $N$ of $P\ga_{\{0\}}P$. By uniqueness, the Fermion complement of $N$ in $P\ga P$ is all of $\bar N\bigvee P\ga_{\{1\}}P$. Thus, in order to shorten the notations, we can suppose that 
$\bar N$ is the Fermion complement of $N$ in $P\ga P$. Thus,
$$
P\ga P=N\bigvee\b^{-1}(\bar N)\sim N\otimes\b^{-1}(\bar N)\,,
$$
where $\b:\tilde N\to\bar N$ is the isomorphism given in \eqref{bbeta}. Define $\ve:=E^{\F\circ\b}_N$ as the Fubini mapping given in \cite{S}, 9.8.4. Let now $x\in N$, $y\in\bar N$. We get
\begin{align*}
&\ve(xy)=\ve(x(y_++y_-))=\ve(xy_+)+\ve[(xV)(Vy_-)]\\
=&\F(y_+)x+\F(\b(Vy_-))xV=\F(y_+)x+\F(y_-)xV\\
=&\F(y_+)x=\F(y_+)x+\F(y_-)x=\F(y)x
\end{align*}
as, being $\F$ even, it is zero on the odd part of $\bar N$.\\
(ii) Take $N_{i}:=P_{i}\car(\ve)P_{i}$, 
$M_{i}:=P_{i}\big(\car(\ve)'\bigwedge\ga\big)P_{i}$, $i=1,2$. 
As $\ve$ is even, we have
\begin{equation*}
\Th(N_{1})=N_{2}\,,\quad
\Th(M_{1})=M_{2}\,,\quad
\Th(P_{1}\ga P_{1})=P_{2}\ga P_{2}\,,
\end{equation*}
and
\begin{align*}
P_{1}\ga P_{1}+P_{2}\ga P_{2}
=&N_{1}\bigvee M_{1}+N_{2}\bigvee M_{2}\\
\sim&N_{1}\otimes M_{1}\oplus N_{2}\otimes M_{2}\,.	
\end{align*}
As $\ve$ is uniquely determined by the restriction on the reduced algebras $\ga_{P_i}$, $i=1,2$, according to \eqref{mmppr} and \eqref{cond}, there exist uniquely determined states $\f_{i}$ on
$M_{i}$, such that
\bes
\ve(x_{1}y_{1}+x_{2}y_{2})=
\f_{1}(y_{1})x_{1}+\f_{2}(y_{2})x_{2}
\ees
whenever $x_{i}\in N_{i}$, $y_{i}\in M_{i}$, $i=1,2$. Thus, it is enough to show that $\f_2=\f_1\circ\Th$.
We compute 
$$
\ve(\Th(x_{1}y_{1}+x_{2}y_{2}))=\f_{1}(\Th(y_{2}))\Th(x_{2})+
\f_{2}(\Th(y_{1}))\Th(x_{1})\,,
$$
and
$$
\Th(\ve(x_{1}y_{1}+x_{2}y_{2}))=\f_{2}(y_{2})\Th(x_{2})+
\f_{1}(y_{1})\Th(x_{1})\,.
$$
Thanks to the $\Th$--equivariance of $\ve$, we conclude that 
$\f_{2}=\f_{1}\circ\Th$ and vice versa.

Finally, if $z\in\ga_{\{1\}}$ is even, then 
$$
P_{i}zP_{i}\equiv P_{i}z\in P_{i}\big(\car(\ve)'
\bigwedge\ga_{[0,1]}\big)P_{i}=M_{i}\,,\quad i=1,2\,.
$$
By the first part, we get
\begin{align*}
&\ve((P_{1}+P_{2})z)=\F(P_{1}zP_{1})P_{1}+\F(\Th(P_{2}zP_{2}))P_{2}\\
=&\F(P_{1}zP_{1})P_{1}+\F(P_{1}zP_{1})P_{2}=\F(P_{1}zP_{1})(P_{1}+P_{2})\,.
\end{align*}
\end{proof}
The previous results relative to the action of the grading automorphism on the centers of the transition expectations allows us to provide the definition of the {\it strongly even} and {\it minimal} Markov states.

We start by noticing that Definition \ref{mrec} can be slightly generalized by simply requiring that the the local algebras $\ga_j$ appearing in \eqref{feal} are full matrix $C^*$--algebras such that the grading automorphism $\Th$ leaves each algebra $\ga_j$ globally 
stable. In this case, the Markov property is still described by the transition expectations $\ve_n$ previously described, and Lemma \ref{orb} still works. 
\begin{defin}
\label{sems}
Let $\f\in\cs(\ga)$ be a Markov state. It is called strongly even (resp. minimal) if the parity automorphism 
$\Th$ acts trivially (resp. transitively) on each $\cz(\car(\ve_{n}))$, 
$\ve_{n}$ being the transition expectations canonically associated to $\f$ through Proposition 
\ref{twst10}. 
\end{defin}
For some interesting applications (see e.g. Corollary \ref{srevreg}), it is enough to consider a Markov state as strongly even if $\Th$ acts trivially on the centers of the transition expectations, infinitely often. Then a Markov state $\f$ will be {\it non strongly even} if there exists $k\in\bn$ such that the action of 
$\Th$ on $\cz(\car(\ve_{n}))$ is nontrivial for each $|n|>k$.

\section{the structure of general Fermi Markov states}

In the present section we investigate the structure of Fermi Markov 
states. We follow Section 3 of \cite{AF1}, where we dealt with the quasi local algebra 
${\displaystyle\ga=\overline{\bigotimes_{j\in I}\bm_{d^j}(\bc)}^{C^*}}$, equipped with the local structure  
${\displaystyle\ga_{\L}=\bigotimes_{j\in\L}\bm_{d^j}(\bc)}$, $\L\subset I$ finite, and trivial $\bz_2$--grading. The forthcoming analysis also represents the extension to the most general Fermion algebra of the results in Section 5 of \cite{AFM}, where only the homogeneous situation 
${\displaystyle\ga:=\overline{\bigvee_{I}\bm_{2}(\bc)}^{C^*}}$, and only the strongly even Markov states were considered. 

The program in \cite{AF1} cannot be directly 
implemented in this situation. In fact 
the parity automorphism $\Th$ does not act trivially on the centres of $\cz(\car(\ve_{j}))$ in general. Thus,  the minimal projections of the centers $\cz(\car(\ve_{j}))$ of 
the ranges 
$\car(\ve_{j})$ does not generate an Abelian algebra. Yet, we are able to decompose 
non homogeneous Markov 
states on the Fermion algebras into minimal ones (cf. Definition \ref{sems}).

Let $\f$ be a Markov state, together with the 
sequence $\{\ve_{j}\}_{j<j_{+}}$ of transition expectations 
canonically associated to $\f$ as previously explained. 
We start by considering the centre $Z_{j}$ of $\car(\ve_{j})$, 
together with the generating family 
$\{P^{j}_{\g_{j}}\}_{\g_{j}\in\G_{j}}$
of minimal projections. Define $\Om_{j}:=\G_{j}/\sim$ where `$\sim$' 
stands for the equivalence relation induced by $\Th$ on the spectrum 
$\G_{j}$ of $Z_{j}$. Let $p_{j}:\G_{j}\to\Om_{j}$ be the 
corresponding canonical projection. Put
\be
\lb{quj}
Q^{j}_{\om_{j}}:=\bigvee_{\g_{j}=p_{j}^{-1}(\{\om_{j}\})}P^{j}_{\g_{j}}\,.
\ee
For $j<j_{+}$, denote $C_{j}\subset Z_{j}$ the subalgebra generated by 
$\{Q^{j}_{\om_{j}}|\om_{j}\in\Om_{j}\}$. Notice that 
$\spec(C_{j})=\Om_{j}$. The $C_{j}$ generate an Abelian 
subalgebra of $\ga$ whose spectrum is precisely
\begin{equation}
\label{compact}
\Om:=\prod_{j\in I}\Om_{j}
\equiv\prod_{j\in I}\sp(C_{j})\,,
\end{equation}
where the product in \eqref{compact} stands for the topological product of the finite sets $\Om_j$.
In order to simplify the notations, we define
$\ve_{j_{+}}:=\id_{\ga_{j_{+}}}$. This means
$\Om_{j_{+}}:=\{j_{+}\}$, $Q_{j_{+}}\equiv P_{j_{+}}:=\idd$, and 
finally, for $N_{j_{+}}$, $\bar N_{j_{+}}$ given in Proposition \ref{part},   
$N_{j_{+}}:=\ga_{\{_{j_{+}}\}}$,  $\bar N_{j_{+}}:=\bc\idd$ with an obvious 
meaning.
Put
$$
B_{j}:=\bigoplus_{\om_{j}\in\Om_{j}}Q^{j}_{\om_{j}}\ga_{\{j\}}Q^{j}_{\om_{j}}\subset\ga_{\{j\}}\,,
$$
and 
\begin{equation*}
\gb:=\overline{\bigvee_{j\in I}B_{j}}\subset\ga\,.
\end{equation*}
The next step is to construct a conditional expectation of $\ga$ onto $\gb$. Thanks to the fact that the $Q^{j}_{\om_{j}}$ are even and thus mutually commuting, we have for each 
$x\in\ga_{[k,l]}$,
\begin{align*}
&\sum_{\om_{k-1},\om_{k},\dots,\om_{l},\om_{l+1}}
\big(Q^{k-1}_{\om_{k-1}}Q^{k}_{\om_{k}}\cdots
Q^{l}_{\om_{l}}Q^{l+1}_{\om_{l+1}}\big)x\big(Q^{k-1}_{\om_{k-1}}Q^{k}_{\om_{k}}\cdots
Q^{l}_{\om_{l}}Q^{l+1}_{\om_{l+1}}\big)\\
=&\sum_{\om_{k-1}}Q^{k-1}_{\om_{k-1}}\sum_{\om_{l+1}}Q^{k}_{\om_{l+1}}
\sum_{\om_{k},\dots,\om_{l}}
\big(Q^{k}_{\om_{k}}\cdots
Q^{l}_{\om_{l}}\big)x\big(Q^{k}_{\om_{k}}\cdots
Q^{l}_{\om_{l}}\big)\\
=&\sum_{\om_{k},\dots,\om_{l}}
\big(Q^{k}_{\om_{k}}\cdots
Q^{l}_{\om_{l}}\big)x\big(Q^{k}_{\om_{k}}\cdots
Q^{l}_{\om_{l}}\big)\,.
\end{align*}
Moreover, if $x=x_{k}x_{k+1}\cdots x_{l}$, then
\begin{align*}
&\sum_{\om_{k},\om_{k+1},\dots,\om_{l},\om_{l}}
\big(Q^{k}_{\om_{k}}Q^{k+1}_{\om_{k+1}}\cdots
Q^{l}_{\om_{l}}\big)x\big(Q^{k}_{\om_{k}}Q^{k+1}_{\om_{k+1}}\cdots
Q^{l}_{\om_{l}}\big)\\
=&\sum_{\om_k,\om_{k+1},\dots,\om_{l}}
\big(Q^{k}_{\om_{k}}x_{k}Q^{k}_{\om_{k}}\big)
\big(Q^{k+1}_{\om_{k+1}}x_{k+1}Q^{k+1}_{\om_{k+1}}\big)
\cdots
\big(Q^{l}_{\om_{l}}x_{l}Q^{l}_{\om_{l}}\big)\,.
\end{align*}
Thus, on the dense subalgebra ${\displaystyle\ca:=\bigcup_{\L\subset I}\ga_{\L}}$, $\L$ finite, we get a norm one projection $E:\ca\to\gb$,
given on the algebraic generators of $\ca$ by
\begin{equation}
\lb{eex}
E\big(x_{j_{1}}\cdots x_{j_{n}}\big)=
\sum_{\om_{j_{1}},\dots,\om_{j_{n}}}
\big(Q^{j_{1}}_{\om_{j_{1}}}x_{j_{1}}Q^{j_{1}}_{\om_{j_{1}}}\big)
\cdots
\big(Q^{j_{n}}_{\om_{j_{n}}}x_{j_{n}}Q^{j_{n}}_{\om_{j_{n}}}\big)
\end{equation}
which uniquely extends to a conditional expectation (denoted again by $E$ by an abuse of notations) $E:\ga\to\gb$ of $\ga$ onto $\gb$.
It is 
also a quite standard fact to see that
\bes
\f=\f\circ E\equiv\f\lceil_{\gb}\,\circ E\,.
\ees
By taking into account the previous considerations we can investigate the structure of Fermi Markov states following the lines in \cite{AF1}. We 
recover the following objects canonically associated to the Markov 
state $\f$ under consideration.
\begin{itemize}
\item[{\bf(a)}] A classical Markov process on the compact space
$\Om$ given in \eqref{compact},
whose law $\m$ is uniquely determined by
the sequences of compatible distributions and transition probabilities
at the place $j$ given respectively by
\begin{align}
\label{isi3}
&\pi^{j}_{\om_{j}}:=
\f(Q^{j}_{\om_{j}})\,,\quad j<j_{+}\\
&\pi^{j}_{\om_{j},\om_{j+1}}:=
\frac{\f(\ve_{j}(Q^{j}_{\om_{j}}Q^{j+1}_{\om_{j+1}}))}
{\f(Q^{j}_{\om_{j}})}\,,\quad j<j_{+}\nn\,.
\end{align}
\item[{\bf(b)}] For each trajectory $\om\equiv(\dots,\om_{j-1},\om_{j},\om_{j+1},\dots)\in\Om$, the $C^{*}$--algebra $\gb^{\om}$ given by
\begin{equation}
\label{disd}
\gb^{\om}:=\overline{\bigvee_{j\in I}
Q^{j}_{\om_{j}}\ga_{\{j\}}Q^{j}_{\om_{j}}}^{C^{*}}\,.
\end{equation}
\end{itemize}
Notice that, in the non trivial cases 
(i.e. when $I$ is infinite), $\gb^{\om}$
cannot be viewed in a canonical way as a subalgebra of $\ga$. Yet, whenever $\L\subset I$ is finite,
$$
\gb^{\om}_{\L}:=\bigvee_{j\in\L}
Q^{j}_{\om_{j}}\ga_{\{j\}}Q^{j}_{\om_{j}}
$$
is a subalgebra of $\ga_{\L}$ with the identity the projection 
${\displaystyle\bigvee_{j\in\L}Q^{j}_{\om_{j}}}$. Namely, 
$\gb^{\om}$ is equipped with a canonical localization 
$\{\gb^{\om}_{\L}\,|\,\L\,\text{finite subset of $I$}\}$, and a 
$\bz_{2}$--grading implemented by the automorphism $\Th^{\om}$ arising from 
the restrictions $\Th\lceil_{\ga_{\L}}$.
\begin{itemize}
\item[{\bf(c)}] A completely positive identity preserving
map $E^{\om}:\ga\to\gb^{\om}$, which is uniquely determined as in \eqref{eex} by
\begin{align}
\label{dise}
x_{j_{1}}x_{j_{2}}\cdots x_{j_{n}}\in\ga\mapsto
&\big(Q^{j_{1}}_{\om_{j_{1}}}x_{j_{1}}Q^{j_{1}}_{\om_{j_{1}}}\big)
\big(Q^{j_{2}}_{\om_{j_{2}}}x_{j_{2}}Q^{j_{2}}_{\om_{j_{2}}}\big)\nn\\
\cdots&
\big(Q^{j_{n}}_{\om_{j_{n}}}x_{j_{n}}Q^{j_{n}}_{\om_{j_{n}}}\big)\,.
\end{align}
\end{itemize}
The above map satisfies $E^{\om}\circ E=E^{\om}$.
\begin{itemize}
\item[{\bf(d)}] A sequence $\{\ce^{\om}_{j}\}_{j\in I}$ of
maps
$$
\ce^{\om}_{j}:\gb^{\om}_{j+1]}\to \gb^{\om}_{j]}
$$
given, for $x\in\gb^{\om}_{j-1]}$, $y\in\gb^{\om}_{[j,j+1]}$ by
\bes
\ce^{\om}_{j}(xy):=\frac{x\ve_{j}(y)}{\pi_{\om_{j},\om_{j+1}}^{j}}\,.
\ees
\end{itemize}
\begin{prop}
\lb{ucx1}
The maps $\ce^{\om}_{j}$ are even conditional expectations.
\end{prop}
\begin{proof}
As for each $k\leq j$, 
$$
\gb^{\om}_{[k,j+1]}=\big(\prod_{l=k}^{j+1}Q^{l}_{\om_{l}}\big)\ga_{[k,j+1]}
\big(\prod_{l=k}^{j+1}Q^{l}_{\om_{l}}\big)\subset\ga_{[k,j+1]}\,,
$$
and
$$
\ce^{\om}_{j}\lceil_{\gb^{\om}_{[k,j+1]}}=\frac{\ce_{j}\lceil_{\gb^{\om}_{[k,j+1]}}}
{\pi_{\om_{j},\om_{j+1}}^{j}}\,.
$$
Thanks to Proposition \ref{twst10}, 
$\ce^{\om}_{j}$ is an even conditional expectation, provided that it 
is identity preserving. This means that we must check 
$\ce^{\om}_{j}(Q^{j}_{\om_{j}}Q^{j+1}_{\om_{j+1}})=Q^{j}_{\om_{j}}$. 
But, we have by \eqref{stpar2} that $\ve_{j}(Q^{j}_{\om_{j}}Q^{j+1}_{\om_{j+1}})=
cQ^{j}_{\om_{j}}$. The proof follows as the number $c$ is precisely 
$\pi_{\om_{j},\om_{j+1}}^{j}$.  
\end{proof}
\begin{itemize}
\item[{\bf(e)}] The state $\psi^{\om}\in\cs(\gb^{\om})$, uniquely 
determined on localized elements by
\begin{equation}
\label{disc}
\psi^{\om}:=\lim_{\stackrel{k\downarrow j_{-}}{l\uparrow j_{+}}}
\frac{\f\big\lceil_{Q^{k}_{\om_{k}}\ga_{\{k\}}Q^{k}_{\om_{k}}}\circ
\ce^{\om}_{k}\circ\cdots\circ\ce^{\om}_{l}}{\pi_{\om_{k}}^{k}}\,.
\end{equation}
\end{itemize}
It is straightforward to check that the state $\psi^{\om}$ is a minimal Markov state on 
$\gb^{\om}$ w.r.t. the conditional expectations 
$$
\widetilde{\ce}^{\om}_{j}:=\ce^{\om}_{j}\circ\ce^{\om}_{j+1}\,.
$$
In addition, the field
$$
\om\in\Om\mapsto\psi^{\om}\circ E^{\om}\in\cs(\ga)
$$
is $\s(\ga^{*},\ga)$--measurable.
\begin{thm}
\label{dis}
Let $\f$ be a Markov state on the Fermion algebra $\ga$
w.r.t. the associated sequence $\{\ce_{j}\}_{j_{-}\leq j<j_{+}}$ of
conditional expectations given in \eqref{cecp}.
Define the compact set $\Om$ by (\ref{compact}), the
probability measure $\m$ on $\Om$  by \eqref{isi3},
the quasi local algebra $\gb^{\om}$ by (\ref{disd}), the
map $E^{\om}$ by (\ref{dise}), the
state $\psi^{\om}$ on $\gb^{\om}$ by (\ref{disc}).

Then $\f$
admits the integral decomposition
\be
\lb{inde}
\f(A)=\int^{\oplus}_{\Om}\psi^{\om}\circ E^{\om}(A)\m(\di\om)\,,\quad A\in\ga\,.
\ee
\end{thm}
\begin{proof}
We outline the proof which is similar to that of Theorem 3.2 of \cite{AF1} after writing down the corresponding objects relative to the Fermi case. 
Consider the Abelian $C^{*}$--subalgebra $\gc$ of $\gb$ given by 
$$
\gc:=\overline{\bigvee_{j\in I}C^{j}}
\sim\overline{\bigotimes_{j\in I}C^{j}}^{C^{*}}\,,
$$ 
together its spectrum $\spec(\gc)=\Om$.
By restricting $\f$ to $\gc$, we obtain a possibly nonhomogeneous Markov random process on $\Om$
with law $\m$ described above.
Let $\pi$ be 
the GNS representation of $\gb$ relative to $\f\lceil_{\gb}$. 
Then 
$L^{\infty}(\Om,\m)\sim\pi(\gc)''\subset\pi(\gb)'\cap\pi(\gb)''$. Thus, we 
have for
$\pi$ the direct integral decomposition
$
{\displaystyle\pi=\int^{\oplus}_{\Om}\pi_{\om}\m(\di\om)}$,
where $\om\mapsto\pi_{\om}$ is a measurable field of representations of
$\gb$. This leads to the direct integral decomposition of $\f_{\lceil\gb}$, and then the decomposition of
$\f\equiv\f\lceil_{\gb}\circ E$ as ${\displaystyle\f=\int_{\Om}\f_{\om}\m(\di\om)}$, see e.g. \cite{T}, Section IV. It is then straightforward to see that 
$$
\f_{\om}(A)=\psi_{\om}(E_{\om}(A))\,,
$$ 
almost everywhere on $\Om$ 
for each $A\in\ga$.
\end{proof}
The constructive part of Proposition \ref{part} allows us to provide the following reconstruction theorem for the class of Fermi Markov states  
considered in the sequel. It parallels the analogous one (cf. \cite{AF1}, Theorem 3.3).

Let $\ga$ be a Fermion algebra. Take for every $j<j_{+}$, a 
$\Th$--invariant commutative subalgebra $Z_{j}$
of $\ga_{\{j\}}$ with spectrum $\G_{j}$ and generators 
$\{P^{j}_{\g_{j}}\}_{\g_{j}\in\G_{j}}$. Put $Z_{j_{+}}:=\bc\idd$.
Let ``$\sim$'' be the equivalence 
relation on the $\G_{j}$ induced by the action of $\Th$, and $p_{j}$ the 
corresponding projection map. Set
$\Om_{j}:=\G_{j}/\sim$, and define $Q^{j}_{\om_{j}}$ as in \eqref{quj}.
Choose a full matrix subalgebra $N^{j}_{\g_{j}}\subset P^{j}_{\g_{j}}\ga_{\{j\}}P^{j}_{\g_{j}}$ 
which is $\Th$--invariant whenever $P^{j}_{\g_{j}}$ 
is a fixed point of $\Th$.\footnote{Notice that $\bar N^{j}_{\g_{j}}$ given in   
Proposition \ref{part} is also left globally invariant under the 
parity.}
Form for $j<j_{+}$, the 
two point even, faithful conditional expectations
\begin{align*}
&\ve^{j}_{\om_{j},\om_{j+1}}:Q^{j}_{\om_{j}}\ga_{\{j\}}Q^{j}_{\om_{j}}
\bigvee Q^{j+1}_{\om_{j+1}}\ga_{\{j+1\}}Q^{j+1}_{\om_{j+1}}\\
&\to\bigoplus_{\g_{j}=p_{j}^{-1}(\{\om_{j}\})}N^{j}_{\g_{j}}
\subset Q^{j}_{\om_{j}}\ga_{\{j\}}Q^{j}_{\om_{j}}\,,\\
&\ve^{j_{+}-1}_{\om_{j_{+}-1},j_{+}}:
Q^{j_{+}-1}_{\om_{j_{+}-1}}\ga_{\{j_{+}-1\}}Q^{j_{+}-1}_{\om_{j_{+}-1}}
\bigvee\ga_{\{j_{+}\}}\\
&\to\bigoplus_{\g_{j_{+}-1}=p_{j}^{-1}(\{\om_{j_{+}-1}\})}N^{j_{+}-1}_{\g_{j_{+}-1}}
\subset Q^{j_{+}-1}_{\om_{j_{+}-1}}\ga_{\{j_{+}-1\}}Q^{j}_{\om_{j_{+}-1}}\,,
\end{align*}
according to Proposition \ref{part}, by taking for the states 
in \eqref{stpar}, \eqref{stpar1}, faithful ones. Define $\gb^{\om}$, $E^{\om}$ 
as in \eqref{disd}, \eqref{dise} respectively. 
For the trajectory $\om=(\dots,\om_{j-1},\om_{j},\om_{j+1},\dots)$, and 
$j<j_{+}$, define the map $\ce_{j}^{\om}$ as
$$
\ce_{j}^{\om}(xy):=x\ve^{j}_{\om_{j},\om_{j+1}}(y)\,,
$$
which is an even faithful conditional expectation according to Proposition \ref{twst10}, and Lemma 
\ref{fcex}.
Take, for $j<j_{+}$, a compatible sequence of even faithful states 
$\f_{j}^{\om}$ on 
$Q^{j}_{\om_{j}}\ga_{\{j\}}Q^{j}_{\om_{j}}$.\footnote{It can be shown 
by a standard compactness property (cf. \cite{AF4}, Proposition 5.1),
that the set of 
sequences of even compatible states $\f_{j}^{\om}$, that is such that 
$\f_{j+1}^{\om}=\f_{j}^{\om}\circ\ce_{j}^{\om}\lceil_{\gb^{\om}_{\{j+1\}}}$,
is nonvoid.} Form the state $\psi^{\om}\in\cs(\gb_{\om})$ as in 
\eqref{disc} by taking as initial distributions the $\f_{j}^{\om}$.
Finally, fix a Markov process on the product space 
${\ds\Om:=\prod_{j\in I}\Om_{j}}$ 
with law $\m$ determined, 
for $\om_{j}\in\Om_{j}$,
$\om_{j+1}\in\Om_{j+1}$, by the marginal
distributions $\pi^{j}_{\om_{j}}>0$, and transition probabilities
$\pi^{j}_{\om_{j},\om_{j+1}}>0$.
\begin{thm}
\label{disbis}
In the above notations, the state $\f$ on $\ga$ given by
$$
\f:=\int_{\Om}\psi_{\om}\circ E_{\om}\m(\di\om)
$$
is a Markov state w.r.t. the sequence $\{\ce_{j}\}_{j<j_{+}}$
of conditional expectations uniquely determined (with the convention $\ga_{\{j_--1\}}=\bc\idd$) for 
$a\in\ga_{j-1]}$, $x\in\ga_{\{j\}}$, $y\in\ga_{\{j+1\}}$ 
by 
\begin{align*}
&\ce_{j}(axy)=a\sum_{\om_{j},\om_{j+1},\om_{j+2}}\pi^{j}_{\om_{j}}
\pi^{j}_{\om_{j},\om_{j+1}}\pi^{j+1}_{\om_{j+1},\om_{j+2}}\\
\times&\ce^{\om}_{j}\big(Q^{j}_{\om_{j}}x\ce^{\om}_{j+1}\big(Q^{j+1}_{\om_{j+1}}y
Q^{j+1}_{\om_{j+1}}Q^{j+2}_{\om_{j+2}}\big)Q^{j}_{\om_{j}}\big)\,,
\qquad j\leq j_{+}-2\,,\\
&\ce_{j}(axy)=a\sum_{\om_{j}}\pi^{j}_{\om_{j}}
\ce^{\om}_{j}\big(Q^{j}_{\om_{j}}xQ^{j}_{\om_{j}}y\big)\,,\qquad\quad\,\,\,\, 
j=j_{+}-1\,.
\end{align*}
\end{thm}
\begin{proof}
A straightforward computation shows that, for all generators of the form 
$x_{k}\cdots x_{l}\in\ga_{[k,l]}$, $\f$ satisfies \eqref{mrec}, for 
the sequence of conditional expectations constructed as above (cf.
Theorem 4.1 of \cite{AF1}). The proof follows as the state $\f$
is locally faithful, by taking into account Lemma \ref{fcex}.
\end{proof}

\section{general properties of Fermi Markov states}

Let $\f\in\cs(\ga)$, and $D_{[k,l]}$ be the adjusted density matrix of the restriction $\f_{[k,l]}$. For 
$k<n<j_{+}$, define 
the unitary $w_{k,n}(t)\in\ga_{[k,n+1]}$ as
\begin{equation*}
w_{k,n}(t):=D_{[k,n+1]}^{it}D_{[k,n]}^{-it}\,,\quad t\in\br\,.
\end{equation*}
The unitaries $\{w_{k,n}(t)\}_{t\in\br}$ give rise to a cocycle 
called {\it transition cocycle} when $\f$ is a Markov state (cf \cite{AFr}).
Denote $S(\,\cdot\,)$ the von Neumann entropy (see e.g. \cite{OP}).

The following theorem collects some properties of the Fermi Markov states, 
which parallel the analogous ones relative to Markov states on tensor product 
algebras (cf. \cite{AF1, AFr, OP}). For the natural applications of the properties described below to the variational principle in quantum statistical mechanics, the reader is referred to \cite{AM2, OP}.
\begin{thm}
\lb{afrl}
Let $\f\in\cs(\ga)$ be a locally faithful even state. Then the following assertions 
are equivalent.
\begin{itemize}
\item[(i)] $\f\in\cs(\ga)$ be a Markov state;
\item[(ii)] for each $t\in\br$ and $k<n<j_{+}$, $w_{k,n}(t)\in\ga_{[n,n+1],+}$.
\end{itemize}
Moreover, if $I=\bz$, $\ga_{\{n\}}=\bm_{2^{d}}(\bc)$ for each $n\in\bz$, and $\f$ is translation invariant,
the previous assertions are also equivalent to 
\begin{itemize}
\item[(iii)] $S(\f_{[0,n+1]})-S(\f_{[0,n]})=S(\f_{[0,1]})-S(\f_{\{0\}})$, $n\geq1$.
\end{itemize}
\end{thm}
\begin{proof}
(i)$\Rightarrow$(ii) Thanks to Lemma 4.1 of \cite{AFr}, 
if $\f$ is a Markov state, then there exists 
an unitary 
${\ds u_{t}\in\ga_{[k,n-1]}'\bigwedge\ga_{[k,n+1]}}$ such that, for 
each $x\in\ga_{[k,n-1]}$,
$$
w_{k,n}(t)xw_{k,n}(t)^{*}\equiv\s^{\f_{[k,n+1]}}_{-t}\big(\s^{\f_{[k,n]}}_{t}(x)\big)
=u_{t}xu_{t}^{*}\equiv x\,,
$$
$\s^{\f}_{t}$ denoting the modular group of a faithful state $\f$ on a von 
Neumann algebra, see e.g. \cite{SZ}.
As $w_{k,n}(t)$ is even, we have
$$
w_{k,n}(t)\in\ga_{[k,n-1]}'\bigwedge\ga_{[k,n+1]}\bigwedge\ga_{+}=\ga_{[n,n+1],+}\,,
$$
see Lemma 11.1 and Theorem 4.17 of \cite{AM2}.

(ii)$\Rightarrow$(i) The Accardi--Cecchini $\f$--expectation 
$E_{k,n}$ of $\f_{[k,n+1]}$ w.r.t. the 
inclusion $\ga_{[k,n]}\subset\ga_{[k,n+1]}$ (cf. \cite{AC}) is written 
as
$$
E_{k,n}(x)=
\ce^{0}_{[k,n]}(w_{k,n}(-i/2)^{*}xw_{k,n}(-i/2))
$$
where $w_{k,n}(-i/2)$ is the analytic continuation at $-i/2$ of 
$w_{k,n}(t)$, and $\ce^{0}_{[k,n]}$ is the conditional expectation of 
$\ga_{[k,n+1]}$ onto 
$\ga_{[k,n]}$ preserving the normalized trace. 
If the $w_{k,n}(t)$ satisfy
all the properties listed above, the Accardi--Cecchini expectation $E_{k,l}$ is a 
$\f_{[k,n+1]}$--preserving quasi conditional expectation
w.r.t. the triplet 
$\ga_{[k,n-1]}\subset\ga_{[k,n]}\subset\ga_{[k,n+1]}$.
By taking for each fixed $n$ the pointwise limit
$$
\ve_{n}:=\lim_{k\downarrow j_{-}}\bigg(\lim_{L}\frac{1}{L}\sum_{l=0}^{L-1}
\big(E_{k,n}\lceil_{\ga_{[n,n+1]}}\big)^{l}\bigg)\,,
$$
we obtain by $\ce_{n}(xy):=x\ve_{n}(y)$, 
$x\in\ga_{[n-1]}$, $y\in\ga_{[n,n+1]}$, a conditional expectation (cf. Proposition \ref{twst10})
fulfilling all the properties listed in Definition \ref{dema}.

(ii)$\iff$(iii) We have 
$$
w_{0,n}(t)=\big[D\f_{[0,n+1]}:D(\f_{[0,n+1]}\circ\ce^{0}_{[0,n]})\big]_{t}\,,
$$
the last being the Connes--Radon--Nikodym cocycle of $\f_{[0,n+1]}$ w.r.t. 
$\f_{[0,n+1]}\circ\ce^{0}_{[0,n]}$ (cf. \cite{SZ}). The assertion follows from the 
fact that (iii) is equivalent to the fact that $\ga_{[n,n+1]}$ is a 
sufficient subalgebra for both the mentioned states. It turns out to 
be  
equivalent to (ii) by translation invariance, see Proposition 11.5 
and Proposition 9.3 of \cite{OP}.  
\end{proof}
\begin{cor}
\lb{kmscor}
Suppose that $j_{-}\in I$. If $\f\in\cs(\ga)$ is a Markov state, 
then its support in $\ga^{**}$ is central. In addition, $\f$ is 
faithful.
\end{cor}
\begin{proof}
By Theorem \ref{afrl}, the pointwise norm limit
$$
\lim_{n\uparrow j_{+}}D_{[j_{-},n]}^{-it}xD_{[j_{-},n]}^{it}
$$
exists as it is asymptotically constant in $n$, on localized elements. 
Thus, it defines a one parameter
group of automorphisms $t\mapsto\s_{t}$ of $\ga$ which admits, by 
construction, $\f$ as a KMS state. This means that 
$\pi_{\f}(\ga)'\xi_{\f}$ is dense in $\ch_{\f}$, 
$(\pi_{\f},\ch_{\f},\xi_{\f})$ being the GNS triplet of $\f$. 
Furthermore, $\f$ is faithful as $\ga$ is a
simple $C^{*}$--algebra, see \cite{BR1}, Proposition 2.6.17.
\end{proof}
\begin{cor}
\label{srevreg}
Suppose that, for each $n\in I$, there exists a $k(n)\in I$ with 
$k(n)\leq n$, such that $\Th$ acts 
trivially on $\cz(\car(\ve_{k(n)}))$, $\ve_{j}$ being the 
transition expectations associated to the Markov state 
$\f$. Then the assertions in Corollary \ref{kmscor} hold true as well.
\end{cor}
\begin{proof}
By regrouping the local algebras, we can suppose that 
there exists a $j_{0}\in I$ such that, for $j<j_{0}$, $\Th$ acts 
trivially on $\cz(\car(\ve_{j}))$.
Consider for $k<j_{0}$, $l>j_{0}$ the local algebras 
$$
\gam_{[k,l]}:=\gar^{c}_{k}\bigvee\ga_{[k+1,l]}\,,
$$
with $\gar^{c}_{k}$ given in \eqref{rcom}.
The last assertion follows as in Corollary 
\ref{kmscor},
by looking at the transition cocycles of $\f$ relative to the 
new localization $\{\gam_{[k,l]}\}_{k<j_{0}<j}$.  
\end{proof}
Let $\f$ be a translation 
invariant locally faithful state on the Fermion 
algebra $\ga\equiv\ga_{\bz}$. The mean entropy $s(\f)$ of $\f$ (see e.g. 
\cite{OP}) 
is defined as
$$
s(\f):=\lim_{n}\frac{1}{n+1}S(\f_{[0,n]})\,,
$$
$S(\f_{[0,n]})$ being the von Neumann entropy of $\f_{[0,n]}$.
\begin{cor}
We have for the translation invariant Markov state $\f$
$$
s(\f)=S(\f_{[0,1]})-S(\f_{\{0\}})\,.
$$
\end{cor}
\begin{proof}
It immediately follows by (iii) of Theorem \ref{afrl}.
\end{proof}

\section{strongly even Markov states}

In the present section we investigate the structure of  strongly even Markov states (cf. Definition 
\ref{sems}). 
By taking into account the structure of the local densities (or equally well the local Hamiltonians by passing to the logaritm) described in \eqref{as}, the strongly even Markov states can be viewed as the Fermi analogue of the Ising type interactions. In addition, they enjoy a kind of 
local entanglement effect, see Section 4 of \cite{FM} for further details.

Notice that the forthcoming analysis extends to the situation when there 
exists a subsequence $\{n_{j}\}\subset I$ such that $\Th$ acts 
trivially on all the $\cz(\car(\ve_{n_{j}}))$.

We start with the following lemma which is known to the experts.
\begin{lem}
\label{mab}
Let $\gc_n\subset\gb_n$, $n\in\bn$, be an increasing sequence of inclusions of unital 
$C^*$--subalgebras of
${\displaystyle\gb:=\overline{\bigcup_{n\in\bn}\gb_n}}$ satisfying $(\gc_k)'\cap\gb_n=\gc_n$, $k\geq n$. Then ${\displaystyle\gc:=\overline{\bigcup_{n\in\bn}\gc_n}}$ is a maximal Abelian $C^*$--subalgebra of $\gb$.
\end{lem}
\begin{proof}
We have for the commutant $\gc'$ in the ambient algebra $\gb$,
\begin{align*}
\gc'=&\overline{\bigcup_{n\in\bn}\big(\gc'\cap\gb_n\big)}
=\overline{\bigcup_{n\in\bn}\bigg(\bigg(\bigcap_{k\in\bn}(\gc_k)'\bigg)\bigcap\gb_n\bigg)}\\
=&\overline{\bigcup_{n\in\bn}\bigg(\bigcap_{k\in\bn}\bigg((\gc_k)'\cap\gb_n\bigg)\bigg)}
=\overline{\bigcup_{n\in\bn}\gc_n}=\gc\,.
\end{align*}
\end{proof}

Let $\om=(\dots,\om_{j-1},\om_{j},\om_{j+1},\dots)\in\Om$ be a 
trajectory. 
Thanks to part (i) of Proposition \ref{part}, 
\begin{align}
\label{disd22}
&\gb^{\om}\equiv
\overline{\big(\bigvee_{j<j_{+}}Q^{j}_{\om_{j}}\ga_{\{j\}}Q^{j}_{\om_{j}}\big)
\bigvee\ga_{\{j_{+}\}}}^{C^{*}}\\
=&\overline{N^{j_{-}}_{\om_{j_{-}}}\bigvee
\big(\bigvee_{j<j_{+}-1}(\bar N^{j}_{\om_{j}}\bigvee
N^{j+1}_{\om_{j+1}})\big)\bigvee(\bar N^{j_{+}-1}_{\om_{j_{+}-1}}
\bigvee \ga_{\{j_{+}\}})}^{C^{*}}\,,
\end{align}
$N^{j}_{\om_{j}}$ $\bar N^{j}_{\om_{j}}$ providing the (Fermi) 
decompositions of the $Q^{j}_{\om_{j}}\ga_{\{j\}}Q^{j}_{\om_{j}}$ described by \eqref{bbar} in 
Proposition \ref{part}. This decomposition is quite similar to the analogous one 
described in Theorem 3.2 of \cite{AF1}, and generalize the situation treated in Section 5 of \cite{AFM}.
\begin{lem}
\lb{note}
Any maximal 
Abelian subalgebra of $\big(\bar N^{j}_{\om_{j}}
\bigvee N^{j+1}_{\om_{j+1}}\big)_{+}$ is maximal Abelian in
$\bar N^{j}_{\om_{j}}\bigvee N^{j+1}_{\om_{j+1}}$ as well.
\end{lem}
\begin{proof}
Let 
$V\in \bar N^{j}_{\om_{j}}\bigvee N^{j+1}_{\om_{j+1}}$ be any 
selfadjoint unitary implementing $\Th$. Then 
\begin{align*}
&\big(\bar N^{j}_{\om_{j}}
\bigvee N^{j+1}_{\om_{j+1}}\big)_{+}\equiv
\big(\bc E_{1}\oplus \bc E_{-1}\big)'\\
=&E_{1}\big(\bar N^{j}_{\om_{j}}\bigvee N^{j+1}_{\om_{j+1}}\big)E_{1}
\oplus E_{-1}\big(\bar N^{j}_{\om_{j}}\bigvee 
N^{j+1}_{\om_{j+1}}\big)E_{-1}\,,
\end{align*}
$V=E_{1}-E_{-1}$ being the resolution of $V$. 
\end{proof}
\begin{lem}
\lb{prle}
The unnormalized trace of 
$$
R:=N^{k}_{\om_{k}}\bigvee\bar N^{k}_{\om_{k}}\bigvee\cdots\bigvee 
N^{l}_{\om_{l}}\bigvee\bar N^{l}_{\om_{l}}
$$
is the product of the unnormalized traces of the $N^{j}_{\om_{j}}$ and
$\bar N^{j}_{\om_{j}}$, $k\leq j\leq l$.
\end{lem}
\begin{proof}
Put ${\displaystyle R=\big(Q^{j}_{\om_{j}}\ga_{\{j\}}Q^{j}_{\om_{j}}\big)
\bigvee\cdots\bigvee
\big(Q^{j+1}_{\om_{j+1}}\ga_{\{j+1\}}Q^{j+1}_{\om_{j+1}}\big)}$.
By the product property of $\tr_{\ga_{[k,l]}}$, 
we get 
$$
\tr{}\!_{R}=\prod_{j=k}^{l}\tr{}\!_{Q^{j}_{\om_{j}}\ga_{\{j\}}Q^{j}_{\om_{j}}}\,.
$$
Thus, we reduce the situation to the algebra 
$N^{j}_{\om_{j}}\bigvee\bar N^{j}_{\om_{j}}\equiv
Q^{j}_{\om_{j}}\ga_{\{j\}}Q^{j}_{\om_{j}}$. Notice that 
$$
Q^{j}_{\om_{j}}\ga_{\{j\}}Q^{j}_{\om_{j}}=
N^{j}_{\om_{j}}\bigvee\tilde N^{j}_{\om_{j}}
\sim N^{j}_{\om_{j}}\otimes\tilde N^{j}_{\om_{j}}\,,
$$
$N^{j}_{\om_{j}}$, $\tilde N^{j}_{\om_{j}}$ are both globally stable 
under the action of $\Th$, 
$\bar N^{j}_{\om_{j}}=\tilde N^{j}_{\om_{j},+}+V\tilde N^{j}_{\om_{j},-}$,
$V$ being any unitary of $N^{j}_{\om_{j}}$ implementing $\Th$ on 
itself, see Proposition \ref{part}.
As the traces are invariant under any automorphism, we get
\begin{align*}
&\tr{}\!_{Q^{j}_{\om_{j}}\ga_{\{j\}}Q^{j}_{\om_{j}}}
=\tr{}\!_{N^{j}_{\om_{j}}}\tr{}\!_{\tilde N^{j}_{\om_{j}}}\\
=&\bigg(\tr{}\!_{N^{j}_{\om_{j},+}}\circ\frac{\id+\Th}{2}\bigg)
\bigg(\tr{}\!_{\tilde N^{j}_{\om_{j},+}}\circ\frac{\id+\Th}{2}\bigg)\\
=&\bigg(\tr{}\!_{N^{j}_{\om_{j},+}}\circ\frac{\id+\Th}{2}\bigg)
\bigg(\tr{}\!_{\bar N^{j}_{\om_{j},+}}\circ\frac{\id+\Th}{2}\bigg)\\
=&\tr{}\!_{N^{j}_{\om_{j}}}\tr{}\!_{\bar N^{j}_{\om_{j}}}\,.
\end{align*}
\end{proof}
Let the initial distributions 
$\eta^{j_{-}}_{\om_{j_{-}}}\in\cs(N^{j_{-}}_{\om_{j_{-}}})$,
the states 
$\eta^{j}_{\om_{j},\om_{j+1}}
\cs(\bar N^{j}_{\om_{j}}\bigvee N^{j+1}_{\om_{j+1}})$ 
be recovered by $\f$ according to
\eqref{stpar}.\footnote{If $j_{-}$ and/or $j_{+}$ do not belong to 
$I$, they do not appear in the formulae, the last having an obvious 
meaning. In addition, as $\Om_{j_{+}}\equiv\{j_{+}\}$, we use the 
symbology $\eta^{j}_{\om_{j},\om_{j+1}}$ also for
the final distributions $\eta^{j_{+}-1}_{\om_{j_{+}-1},j_{+}}$.}
Consider the even densities $T^{(j)}_{\om_{j}}$, $\hat T^{(j)}_{\om_{j}}$,
$T^{(j)}_{\om_{j},\om_{j+1}}$ localized in $N^{j}_{\om_{j}}$,
$\bar N^{j}_{\om_{j}}$, $\bar N^{j}_{\om_{j}}\bigvee N^{j+1}_{\om_{j+1}}$, 
and associated to $\eta^{j_{-}}_{\om_{j_{-}}}$ or
$\eta^{j}_{\om_{j},\om_{j+1}}\lceil_{N^{j+1}_{\om_{j+1}}}$, 
$\eta^{j}_{\om_{j},\om_{j+1}}\lceil_{\bar N^{j}_{\om_{j}}}$,
$\eta^{j}_{\om_{j},\om_{j+1}}$
respectively.
\begin{prop}
\lb{prle1}
The states $\eta^{j_{-}}_{\om_{j_{-}}}$, $\eta^{j}_{\om_{j},\om_{j+1}}$
uniquely define a product state 
on $\gb^{\om}$, coinciding with $\psi^{\om}$ in 
\eqref{inde}, which is symbolically written as
\begin{equation}
\label{disc22}
\psi^{\om}=\eta^{j_{-}}_{\om_{j_{-}}}\prod_{j\leq j_{+}-1}
\eta^{j}_{\om_{j},\om_{j+1}}\,.
\end{equation}
\end{prop}
\begin{proof}
Consider on $\gb^{\om}$ the localization 
$$
\gb^{\om}=\overline{\bigvee_{j\in I}\gn^j_{\om_j}}
$$
suggested by \eqref{disd22}. Here, $\gn^{j_-}_{\om_{j_-}}:=N^{j_-}_{\om_{j_-}}$, 
$\gn^j_{\om_j}:=\bar N^{j-1}_{\om_{j-1}}\bigvee
N^{j+1}_{\om_{j}}$, $j_-<j<j_+$, and finally 
$\gn^{j_+}_{\om_{j_+}}:=\bar N^{j_{+}-1}_{\om_{j_{+}-1}}
\bigvee \ga_{\{j_{+}\}}$.
As the above densities commute each other, for each $k<l$, the product of local densities
$$
T^{(k-1)}_{\om_{k-1},\om_{k}}
\times\cdots\times
T^{(l-1)}_{\om_{l-1},\om_{l}}
$$
is a well defined positive even operator on $\bigvee_{k\leq j\leq l}\gn^j_{\om_j}$
which by Lemma \ref{prle}, is the density of $\psi^{\om}\lceil_{\bigvee_{k\leq j\leq l}\gn^j_{\om_j}}$
w.r.t. the unnormalized 
trace of $\bigvee_{k\leq j\leq l}\gn^j_{\om_j}$. As explained in Section 2.3, this means that $\psi^{\om}$ is the product states of $\eta^{j_{-}}_{\om_{j_{-}}}$ with the $\eta^{j}_{\om_{j},\om_{j+1}}$ as explained in \eqref{disc22} (see Theorem 11.2 of \cite{AM2} for a similar situation).
\end{proof}
As all the states appearing in \eqref{disc22} are even, we can 
explicitely write the local densities associated to the 
strongly even Markov state. Namely,
consider the Radon--Nikodym derivatives (i.e. the densities) 
$T_{\ga_{[k,l]}}$ w.r.t. the unnormalized 
trace of $\ga_{[k,l]}$,
$$
\f_{[k,l]}=\tr{}\!_{\ga_{[k,l]}}(T_{\ga_{[k,l]}}\,\cdot\,)\,.
$$
Then $T_{\ga_{[k,l]}}$ has the nice decomposition
\begin{equation}
\label{as}
T_{\ga_{[k,l]}}=\bigoplus_{\om_{k},\dots,\om_{l}}
T^{(k)}_{\om_{k}}T^{(k)}_{\om_{k},\om_{k+1}}
\times\cdots\times
T^{(l-1)}_{\om_{l-1},\om_{l}}\hat T^{(l)}_{\om_{l}}\,.
\end{equation}
By Corollary \ref{kmscor}, any strongly even Markov state is a KMS state 
for the one parameter group of automomorphisms $\s_{t}$ given, for 
$x\in\ga$, by
$$
\s_{t}(x):=\lim_{\stackrel{k\downarrow j_{-}}{l\uparrow j_{+}}}
T_{\ga_{[k,l]}}^{-it}xT_{\ga_{[k,l]}}^{it}\,.
$$
In addition, each strongly even Markov state is faithful.

We now show that each strongly even Markov state  
is a lifting of a classical Markov processes. This result parallels the analogous one relative to the tensor product algebra, obtained first in \cite{GZ} for some particular cases, and then in \cite{FM} for the general situation. Such property was called {\it diagonalizability} in \cite{GZ}. After adapting the situation relative to the tensor product case to the strongly even Fermi Markov states, we can follow the same line of the proof of Theorem 3.2 of \cite{FM}.

We start by defining increasing 
subalgebras of the Fermion algebra $\ga$ equipped with a natural 
local structure inherited from that of the original algebra.
Let $\gar_{j}:=\car(\ve_{j})$, with relative commutant 
\be
\lb{rcom}
\gar_{j}^{c}:=\gar_{j}'\bigwedge\ga_{\{j\}}\,.
\ee
Define
\begin{align}
\lb{ntppp}
&\gn_{\{k\}}:=\cz(\gar_{k})\,,\quad\gn_{[k,k+1]}:=\gar_{k}^{c}\bigvee\gar_{k+1}\,,\\
&\gn_{[k,l]}:=\gar_{k}^{c}\bigvee\ga_{[k+1,l-1]}\bigvee\gar_{l}\,,\quad k<l+1\,.\nn
\end{align}
Thanks to Lemma \ref{note}, for each $k\leq j<l$ and $\om_{j}\in\Om_{j}$, we can choose a even
maximal
Abelian subalgebra $D^{j}_{\om_{j},\om_{j+1}}$ of
$\bar N^{j}_{\om_{j}}\bigvee N^{j+1}_{\om_{j+1}}$ containing
$T^{(j)}_{\om_{j},\om_{j+1}}$.
Put
\begin{align}
\lb{ntp}
&\gd_{\{k\}}:=\gn_{\{k\}}\equiv\cz(\gar_{k})\,,\nn\\
&\gd_{[k,l]}:=\bigoplus_{\om_{k},\dots,\om_{l}}
\big(D^{k}_{\om_{k},\om_{k+1}}\bigvee\cdots\bigvee
D^{l-1}_{\om_{l-1},\om_{l}}\big)\,,\quad k<l\,,\\
&\gd:=\overline{\big(\bigcup_{[k,l]\subset I}\gd_{[k,l]}\big)}\,.\nn
\end{align}
\begin{thm}
\label{dmain} 
Let $\f\in\cs(\ga)$ be a strongly even Markov state. Then there 
exists an even maximal Abelian $C^*$--subalgebra $\gd\subset\ga$, and a
conditional expectation $\ge:\ga\to\gd$ such that
$\f=\f\lceil_{\gd}\circ\ge$. 
In addition, the measure $\m$  on $\spec(\gd)$
associated to $\f\lceil_{\gd}$ is a Markov measure w.r.t. the natural 
localization of $\gd$ given in \eqref{ntp}. 
\end{thm}
\begin{proof}
Let $[m_k,n_k]$ be an increasing sequence of intervals such that $[m_k,n_k]\uparrow I$. Then 
$$
\ga=\overline{\big(\lim_{\stackrel{\longrightarrow}
{[m_k,n_k]\uparrow I}}\gn_{[m_k,n_k]}\big)}^{\,C^{*}}\,.
$$
As $\gd_{[m,n]}$ is an even maximal Abelian subalgebra of $\gn_{[m,n]}$, the increasing sequence
$\gd_{[m_k,n_k]}\subset\gn_{[m_k,n_k]}$ satisfies the hypotheses of Lemma \ref{mab}.
Thus, $\gd$ is a even maximal Abelian $C^{*}$--subalgebra of $\ga$. 
According to \eqref{as}, we have 
$$
T_{\gn_{[m,n]}}=\bigoplus_{\om_{m},\dots,\om_{n}}
T^{(m)}_{\om_{m},\om_{m+1}}
\times\cdots\times
T^{(n-1)}_{\om_{n-1},\om_{n}}\,,
$$
that is, $\{T_{\gn_{[m,n]}}\}_{m<n}\subset\gd$. 
Let $E^{0}_{m,n}:\gn_{[m,n]}\to\gd_{[m,n]}$ be the canonical conditional 
expectation
of $\gn_{[m,n]}$ onto the maximal abelian subalgebra
$\gd_{[m,n]}$ (cf. \cite{FM}, Footnote 4). We have 
\begin{align}
\label{iss4}
\f\lceil_{\gn_{[m,n]}}\equiv&\tr{}\!_{\gn_{[m,n]}}
\big(T_{\gn_{[m,n]}}\,\cdot\,\big)
=\tr{}\!_{\gn_{[m,n]}}\big(E^{0}_{m,n}(T_{\gn_{[m,n]}}\,\cdot\,)\big)\nn\\
=&\tr{}\!_{\gn_{[m,n]}}\big(T_{\gn_{[m,n]}}E^{0}_{m,n}(\,\cdot\,)\big)
\equiv\f\lceil_{\gn_{[m,n]}}\circ E^{0}_{m,n}\,.
\end{align}
As the sequence $\{E^{0}_{m,n}\}_{m<n}$ is projective, the direct limit 
${\ds\lim_{\stackrel{\longrightarrow}{[m,n]\uparrow I}}E^{0}_{m,n}}$ 
uniquely defines a conditional expectation $\ge:\ga\to\gd$ 
fulfilling by \eqref{iss4}, $\f=\f\lceil_{\gd}\circ\ge$.
The measure $\m$ on $\spec(\gd)$
associated to $\f\lceil_{\gd}$ is a Markov measure w.r.t. the natural 
localization of $\gd$ previously described.
This follows as in Section 6 of \cite{FM}, 
after noticing that $D^{m}_{\om_{m},\om_{m+1}}\bigvee\cdots\bigvee
D^{n-1}_{\om_{n-1},\om_{n}}$
in \eqref{ntp} generates a tensor product, and the restriction 
$\f\lceil_{D^{m}_{\om_{m},\om_{m+1}}\bigvee\cdots\bigvee
D^{n-1}_{\om_{n-1},\om_{n}}}$ defines a product measure on 
$\spec\big(D^{m}_{\om_{m},\om_{m+1}}\big)\times\cdots\times
\spec\big(D^{n}_{\om_{n},\om_{n+1}}\big)$.
\end{proof}
Now we pass to the dynamical entropy $h_{\f}(\a)$ 
w.r.t. the right shift $\a$ for translation invariant 
strongly even Markov states. The reader is referred to \cite{CP, CNT, 
OP} for the definition and technical details on the dynamical entropy.

The definition of the dynamical entropy $h_{\f}(\a)$ is based on the 
multiple subalgebra entropy $H_{\f}(N_{1},\cdots,N_{k})$, with 
$N_{1},\cdots,N_{k}\subset M$. We start by pointing out that, if the 
subalgebras $N_{1},\cdots,N_{k}$ are the range of $\f$--preserving 
conditional expectations and are contained in different factors of a 
tensor product algebra, then
\be
\lb{cnten}
H_{\f}(N_{1},\cdots,N_{k})=S\big(\f\lceil_{N}\big)\,,
\ee
with $N:=N_{1}\bigvee\cdots\bigvee N_{k}$.\footnote{Fix a faithful 
trace on $M$. Let $T_{1},\cdots,T_{k}$, $T$ be the corresponding 
densities of $N_{1},\cdots,N_{k}$, $M$ respectively. Choose maximal 
Abelian subalgebras $A_{j}$ of $N_{j}$ containing $T_{j}$, 
$j=1,\dots,k$. As the $N_{j}$ are expected, we have for $a\in A_{j}$,
$$
T^{-it}aT^{it}=T_{j}^{-it}aT_{j}^{it}=a\,,
$$
that is $A_{j}\subset M_{\f}$, $M_{\f}$ being the centralizer of the 
faithful state $\f$.
As the $A_{j}$ are contained in different factors of a tensor product, 
$A_{1}\bigvee\cdots\bigvee A_{k}$ is maximal Abelian 
in $N$. Thus, \eqref{cnten} follows by Corollary VIII.8 of \cite{CNT}.} 
\begin{thm}
\lb{cnten1}
Let $\f\in\cs(\ga)$ be a translation invariant strongly even Markov 
state. Then $h_{\f}(\a)=s(\f)$.
\end{thm}
\begin{proof}
The proof follows the same lines of the tensor product case. We keep into account some boundary effects which cannot be neglected in proving the result.
Fix $n$, and consider 
$\gn_{[0,n+1]}$ given in \eqref{ntppp}. We have 
$\ga_{[1,n]}\subset\gn_{[0,n]}\subset\ga_{[0,n+1]}$, and $\gn_{[0,n]}$ 
is expected. We compute,
\begin{align*}
H(k):=&H_{\f}(\gn_{[0,n]},\a(\gn_{[0,n]}),\cdots,\a^{k(n+2)}(\gn_{[0,n]}))\\
\geq&H_{\f}(\gn_{[0,n]},\a^{n+2}(\gn_{[0,n]}),\cdots,\a^{k(n+2)}(\gn_{[0,n]}))\\
\geq&H_{\f}(\gn_{[0,n],+},\a^{n+2}(\gn_{[0,n],+}),\cdots,\a^{k(n+2)}(\gn_{[0,n],+}))\,,
\end{align*}
Now, 
$\gn_{[0,n],+},\a^{n+2}(\gn_{[0,n],+},),\dots,\a^{k(n+2)}(\gn_{[0,n],+})$ 
are all expected, and generate a tensor product. Then
$$
H(k)\geq 
S\big(\f\lceil_{M_{k}}\big)
=-S\big(\f\lceil_{M_{k}},\t\lceil_{M_{k}}\big)
+k\ln d\,.
$$ 
Here, $M_{k}:=\gn_{[0,n],+}\bigvee\a^{n+2}(\gn_{[0,n],+})\bigvee
\dots\bigvee\a^{k(n+2)}(\gn_{[0,n],+})$,
$d$ is the tracial dimension of $\gn_{[0,n],+}$, $\t$ the normalized 
trace on $\ga$, and finally $S(\,\cdot\,,\,\cdot\,)$ the relative 
entropy (see e.g. \cite{OP}). As $\ga_{[1,m],+}\subset\gn_{[0,m],+}$, and the tracial 
dimension of $\ga_{[1,m],+}$ coincides with that of 
$\ga_{[1,m]}$ (cf. Lemma \ref{note}),
we obtain by the monotonicity of the relative entropy, 
\begin{align*}
H(k)\geq-&S\big(\f\lceil_{\ga_{[1,(n+2)(k+1)]}},\t\lceil_{\ga_{[1,(n+2)(k+1)]}}\big)
+k\ln d\\
=&S\big(\f\lceil_{\ga_{[1,(n+2)(k+1)]}}\big)+[kn-(n+2)(k+1)]\ln l\,,
\end{align*}
$l$ being the tracial dimension of $\ga_{\{0\}}$. Finally, we get
\begin{align*}
&h_{\f}(\a)\geq\lim_{k}\frac{H(k)}{(n+2)k}
\geq\lim_{k}\bigg[\frac{k+1}{k}s(\f)\\
+&\frac{kn-(n+2)(k+1)}{(n+2)k}\ln l\bigg]=s(\f)-\frac{2\ln l}{n+2}\,.
\end{align*}
Since $h_{\f}(\a)\leq s(\f)$ and $n$ is arbitrary, the assertion 
follows.
\end{proof}

\section{examples of translation invariant Fermi Markov states}
\label{exxx}

In the present section we exhibit some examples of Fermi Markov states. 
We restrict the matter to the translation invariant 
situation. The non homogeneous cases can be analogously treated.
The present construction furnishes the direct application of Theorem  \ref{disbis}, or equally well Proposition \ref{part}.
Thanks to the translation invariance, it is 
enough to
construct a two
point even transition expectation 
$\ve:\ga_{[0,1]}\to\ga_{\{0\}}$, and compute the stationary 
even distributions by solving
$\r=\r\circ\ve\circ\a\lceil_{\ga_{\{0\}}}$, 
$\r$ running into the even states of $\ga_{\{0\}}$. A translation 
invariant Markov state $\f$ is then recovered by the marginals
\be
\lb{margi}
\f(x_{k}\cdots x_{l})=\r(\ve_{k}(x_{k}\ve_{k+1}(x_{k+1}\cdots\ve_{l-1}(x_{l-1}
\ve_{l}(x_{l}))\cdots)))\,.
\ee

\subsection{Case 1:} $\ga_{\{n\}}\sim\bm_{2}(\bc)$, $\cz(\car(\ve))\sim\bc^{2}$, 
$\#$ of orbits of $\Th\lceil_{\cz(\car(\ve))}=1$.

We start with the pivotal example in Subsection 6.4 of \cite{AFM} by showing that it provides examples of Fermi Markov states which are entangled. Define, for a fixed $\chi$ in the unit circle $\bt$,
\begin{equation*}
q_{\chi}:=\frac{1}{2}\big(\idd+\chi 
a_{0}+\overline{\chi}a_{0}^{+}\big)\,.
\end{equation*}
Choose a faithful state $\eta\in\cs(q_{\chi}\ga_{[0,1]}q_{\chi})$. Put
\begin{equation}
\label{chi1}
\ve(x)=\eta(q_{\chi}xq_{\chi})q_{\chi}
+\eta(q_{\chi}\Th(x)q_{\chi})q_{-\chi}\,,\quad x\in\ga_{[0,1]}\,.
\end{equation}
With $\t$ the normalized trace on $\bm_{2}(\bc)$, 
$\ve_{n}:=\ve\circ\a^{-n}$, and $x_{k}\in\ga_{\{k\}},\dots,\\
x_{l}\in\ga_{\{l\}}$, the marginals \eqref{margi} with $\r=\t$,
uniquely determine a translation invariant locally faithful Markov state 
$\f$ on the Fermion algebra
$\ga:=\ga_{\bz}$ satisfying the required properties. 
Thanks to shift invariance, 
it suffices to consider $x\in\ga_{\{0\}}$, $y\in\ga_{\{1\}}$.

Let $\xi_{\chi}$, $\xi^{\perp}_{\chi}$ be the (uniquely determined up to a phase) eigenvectors of 
$q_{\chi}$, $q_{-\chi}=q^{\perp}_{\chi}$ acting on $\bc^2$, corresponding to the eigenvalues 1, respectively. Put 
$$
V:=\langle\,{\bf\cdot}\,,\xi_{\chi}\rangle\xi^{\perp}_{\chi}\,.
$$ 
As $V\in\bm_2(\bc)=\ga_{\{0\}}\subset\ga_{[0,1]}$, $V$ is also in $\ga_{[0,1]}$. Put 
$\d:=\eta(V(\chi a_{1}+\overline{\chi}a^{+}_{1})q_{\chi}))$.
We have
$$
\f(xy)=\langle x_{+}\xi_{\chi},\xi_{\chi}\rangle\langle\a^{-1}(y_{+})\xi_{\chi},\xi_{\chi}\rangle
+\d\langle x_{-}\xi_{\chi},\xi^{\perp}_{\chi}\rangle\langle\a^{-1}(y_{-})\xi_{\chi},\xi_{\chi}\rangle\,.
$$
Now we show that there exists a faithful state 
$\eta$ as above, such 
that $\eta(X)\neq0$, where  
\begin{equation}
\label{cchi1}
X:=V\a(q_{\chi,-})q_{\chi}\equiv\frac{1}{2}V(\chi a_{1}+
\overline{\chi}a^{+}_{1})q_{\chi})\neq0\,.
\end{equation}
Pick a functional which is different from zero on 
$X$, hence a state $\eta_{0}$ which is nonnull on $X$. Let 
$p\in q_{\chi}\ga_{[0,1]}q_{\chi}$ be the support of $\eta_{0}$. Choose a state 
$\eta_{1}$ with support $q_{\chi}-p$. Then
$\eta:=\b\eta_{0}+(1-\b)\eta_{1}$ is a faithful state on
$q_{\chi}\ga_{[0,1]}q_{\chi}$ which is 
nonnull on $X$ for an appropriate choice of 
$\b\in [0,1]$.\footnote{The last claim easily follows as 
$\eta(X)=0$ means $\eta_{0}(X)\neq\eta_{1}(X)$, and
${\ds\b=\frac{\eta_{1}(X)}{\eta_{1}(X)-\eta_{0}(X)}}$.} We then have the following
\begin{prop}
\label{eeeent}
Let $\L_1,\L_2\subset\bz$ such that $\L_1\bigcap\L_2=\emptyset$, $\L_1\bigcup\L_2=\bz$. Suppose that $\eta(X)\neq0$, where $\eta$ is the state in \eqref{chi1} and $X$ is given in \eqref{cchi1}. Then the state $\f$ described above is entangled w.r.t. the decomposition 
$\ga=\overline{\ga_{\L_1}\bigvee\ga_{\L_2}}$.
\end{prop}
\begin{proof}
Let $\ga_{\{n\}}\subset\ga_{\L_1}$, $\ga_{\{n+1\}}\subset\ga_{\L_2}$ for some $n\in\bz$ (which is always the case after a possible renumbering of $\L_1$, $\L_2$). Under the above assumption, $\f(x_{-}y_{-})$ cannot be identically zero for each $x\in\ga_{\L_1}$, $y\in\ga_{\L_2}$ due to the shift invariance. The proof now follows by applying the Moriya criterion established in Proposition 1 of \cite{M}.
\end{proof}

By extending the previous computations to more general cases, it is then possible to construct many examples of entangled translation invariant Fermi Markov states for the situation when 
$\ga_{\{0\}}=\bm_{2^d}$, $d>1$. We are going to describe a sample of pivotal examples.

We now consider the successive step 
$\ga_{\{k\}}\sim\bm_{4}(\bc)$.  
We exhibit examples for each possible structure of the Abelian 
algebra $\cz(\car(\ve))$, and for the action of $\Th$ on it. Let 
$\{a_{i},a^{+}_{i}\,|\,i=1,2\}$ be the creators and annihilators 
generating $\ga_{\{0\}}$. Consider the system $\{e_{kl}(j)\,|\,j,k,l=1,2\}$
of commuting 
$2\times2$ matrix
units obtained via the Jordan--Klein--Wigner 
transformation \eqref{kw}. Putting 
$e_{(i,j)(k,l)}:=e_{ik}(1)e_{jl}(2)$, we obtain a system of matrix units 
for $\ga_{\{0\}}$ which realizes the isomorphism 
$\ga_{\{0\}}\sim\bm_{2}(\bc)\otimes\bm_{2}(\bc)$.

\subsection{Case 2:} $\cz(\car(\ve))\sim\bc^{4}$, 
$\#$ of orbits of $\Th\lceil_{\cz(\car(\ve))}=4$.

Choose $\{e_{(i,j)(i,j)}\,|\,i,j=1,2\}$ as the generators of 
$\cz(\car(\ve))$. In this situation, there exist even states $\f_{ij}$, 
$i,j=1,2$ on $\ga_{\{1\}}$ such that for $x\in\ga_{\{0\}}$, 
$y\in\ga_{\{1\}}$,
$$
\ve(xy)=\sum_{i,j=1}^{2}\tr(xe_{(i,j)(i,j)})\f_{ij}(y)e_{(i,j)(i,j)}\,.
$$
This is nothing but the example in Subsection 6.2 of \cite{AFM}. Thus, 
$\f$ is strongly clustering w.r.t. the shift on the chain, and the 
von Neumann algebra $\pi_{\f}(\ga)''$ generated by the GNS 
representation $\pi_{\f}$ of $\f$ is a type $\ty{III_{\l}}$ factor for 
some $\l\in (0,1]$, see \cite{FM}.

\subsection{Case 3:} $\cz(\car(\ve))\sim\bc^{4}$, 
$\#$ of orbits of $\Th\lceil_{\cz(\car(\ve))}=3$.

For a fixed $\chi$ in the unit circle $\bt$, define 
$$
Q_{\chi}:=\frac{1}{2}\big(\idd+\chi 
a_{2}+\overline{\chi}a_{2}^{+}\big)\,.
$$
Choose $\{e_{(1,j)(1,j)},e_{22}(1)Q_{\pm\chi}\,|\,j=1,2\}$ as the generators of 
$\cz(\car(\ve))$. In this situation, there exist even states $\f_{j}$, 
$j=1,2$ on $\ga_{\{1\}}$, and a state $\f$ on 
$e_{22}(1)Q_{\chi}\ga_{[0,1]}e_{22}(1)Q_{\chi}$ such that, 
for $x\in\ga_{\{0\}}$, $y\in\ga_{\{1\}}$,
\begin{align*}
\ve(xy)=&\sum_{j=1}^{2}\tr(xe_{(1,j)(1,j)})\f_{j}(y)e_{(1,j)(1,j)}\\
+&\f(e_{22}(1)Q_{\chi}xye_{22}(1)Q_{\chi})
e_{22}(1)Q_{\chi}\\
+&\f(e_{22}(1)Q_{\chi}\Th(xy)e_{22}(1)Q_{\chi})
e_{22}(1)Q_{-\chi}\,.
\end{align*}

\subsection{Case 4:} $\cz(\car(\ve))\sim\bc^{4}$, 
$\#$ of orbits of $\Th\lceil_{\cz(\car(\ve))}=2$.

First choose
$\{e_{ii}(1)Q_{\pm\chi}\,|\,i=1,2\}$ as the generators of 
$\cz(\car(\ve))$. In this situation, there exist states $\f_{i}$, 
on $e_{ii}(1)Q_{\chi}\ga_{[0,1]}e_{ii}(1)Q_{\chi}$, $i=1,2$ 
such that, for $x\in\ga_{[0,1]}$,
\begin{align*}
\ve(x)=&\sum_{i=1}^{2}\big(\f_{i}(e_{ii}(1)Q_{\chi}xe_{ii}(1)Q_{\chi})
e_{ii}(1)Q_{\chi}\\
+&\f_{i}(e_{ii}(1)Q_{\chi}\Th(x)e_{ii}(1)Q_{\chi})
e_{ii}(1)Q_{-\chi}\big)\,.
\end{align*}
Next, for fixed $(\chi,\eta)\in\bt^{2}$, define with 
$V:=a_{1}^{+}a_{1}-a_{1}a_{1}^{+}$,
$$
P_{\chi,\eta}:=\frac{1}{4}\big(\idd+\chi 
a_{1}+\overline{\chi}a_{1}^{+}\big)\big(\idd+\eta V 
a_{2}+\overline{\eta}Va_{2}^{+}\big)\,.
$$
Choose $\{P_{\pm\chi,\pm\eta}\}$ as the generators of 
$\cz(\car(\ve))$. In this situation, there exist states $\f_{\pm}$ 
on $P_{\pm\chi,\eta}\ga_{[0,1]}P_{\pm\chi,\eta}$ respectively, 
such that for $x\in\ga_{[0,1]}$,
\begin{align*}
\ve(x)=&\f_{+}(P_{\chi,\eta}xP_{\chi,\eta})
P_{\chi,\eta}+\f_{+}(P_{\chi,\eta}\Th(x)P_{\chi,\eta})
P_{-\chi,-\eta}\\
+&\f_{-}(P_{-\chi,\eta}xP_{-\chi,\eta})
P_{-\chi,\eta}+\f_{-}(P_{-\chi,\eta}\Th(x)P_{-\chi,\eta})
P_{\chi,-\eta}\,.
\end{align*}

\subsection{Case 5:} $\cz(\car(\ve))\sim\bc^{3}$, 
$\#$ of orbits of $\Th\lceil_{\cz(\car(\ve))}=3$.

First choose
$\{e_{11}(1)e_{jj}(2)\,,e_{22}(1)\,|\,j=1,2\}$ as the generators of 
$\cz(\car(\ve))$. We have two possibilities. Namely, there exist even states 
$\f_{j}$, on $\ga_{\{1\}}$, $i=1,2$, and an even state $\f$ either on 
$\ga_{\{1\}}$, or on $(e_{22}(1)\ga_{\{0\}}e_{22}(1))\bigvee\ga_{\{1\}}$
such that, for $x\in\ga_{\{0\}}$, $y\in\ga_{\{1\}}$,
$$
\ve(xy)=\sum_{j=1}^{2}\tr(xe_{11}(1)e_{jj}(2))\f_{j}(y)e_{11}(1)e_{jj}(2)
+\f(y)e_{22}(1)xe_{22}(1)\,,
$$
respectively
$$
\ve(xy)=\sum_{j=1}^{2}\tr(xe_{11}(1)e_{jj}(2))\f_{j}(y)e_{11}(1)e_{jj}(2)
+\f(e_{22}(1)xe_{22}(1)y)e_{22}(1)\,.
$$
Next, put $P:=e_{(1,2)(1,2)}+e_{(2,1)(2,1)}$ and choose
$\{e_{(i,i)(i,i)}\,,P\,|\,i=1,2\}$ as the generators of 
$\cz(\car(\ve))$. Again, we have two possibilities. Namely, there exist even states 
$\f_{j}$, on $\ga_{\{1\}}$, $i=1,2$, and an even state $\f$ either on 
$\ga_{\{1\}}$, or on $(P\ga_{\{0\}}P)\bigvee\ga_{\{1\}}$
such that, for $x\in\ga_{\{0\}}$, $y\in\ga_{\{1\}}$,
$$
\ve(xy)=\sum_{j=1}^{2}\tr(xe_{11}(1)e_{jj}(2))\f_{j}(y)e_{11}(1)e_{jj}(2)
+\f(y)PxP\,,
$$
respectively
$$
\ve(xy)=\sum_{j=1}^{2}\tr(xe_{11}(1)e_{jj}(2))\f_{j}(y)e_{11}(1)e_{jj}(2)
+\f(PxPy)P\,.
$$
Notice that the last possibilities correspond to nontrivial cases 
with $\car(\ve)\subset\ga_{+}$.

\subsection{Case 6:} $\cz(\car(\ve))\sim\bc^{3}$, 
$\#$ of orbits of $\Th\lceil_{\cz(\car(\ve))}=2$.

For $\chi\in\bt$, choose $\{e_{11}(1)Q_{\pm\chi}\,,e_{22}(1)\}$ as the generators of 
$\cz(\car(\ve))$. We have two possibilities. Namely, choose a 
state $\f$ on $e_{11}(1)Q_{\chi}\ga_{[0,1]}e_{11}(1)Q_{\chi}$, and a 
even state $\psi$ either on $\ga_{\{1\}}$, or on $(e_{22}(1)\ga_{\{0\}}e_{22}(1))
\bigvee\ga_{\{1\}}$
such that, for $x\in\ga_{\{0\}}$, $y\in\ga_{\{1\}}$,
\begin{align*}
\ve(xy)=&\f(e_{11}(1)Q_{\chi}xye_{11}(1)Q_{\chi})e_{11}(1)Q_{\chi}\\
+&\f(e_{11}(1)Q_{\chi}\Th(xy)e_{11}(1)Q_{\chi})e_{11}(1)Q_{-\chi}\\
+&\psi(y)e_{22}(1)xe_{22}(1)\,,
\end{align*}
respectively
\begin{align*}
\ve(xy)=&\f(e_{11}(1)Q_{\chi}xye_{11}(1)Q_{\chi})e_{11}(1)Q_{\chi}\\
+&\f(e_{11}(1)Q_{\chi}\Th(xy)e_{11}(1)Q_{\chi})e_{11}(1)Q_{-\chi}\\
+&\psi(e_{22}(1)xe_{22}(1)y)e_{22}(1)\,.
\end{align*}

\subsection{Case 7:} $\cz(\car(\ve))\sim\bc^{2}$, 
$\#$ of orbits of $\Th\lceil_{\cz(\car(\ve))}=2$.

We treat only the following cases, the remaining ones follow 
analogously. Choose $p=e_{(1,1)(1,1)}$, $p^{\perp}$ as the generators 
of $\cz(\car(\ve))$. We have two possibilities. Namely, there exists a 
even state $\f$ on $\ga_{\{1\}}$, and a 
even state $\psi$ either on $\ga_{\{1\}}$, or on $(p^{\perp}\ga_{\{0\}}p^{\perp})
\bigvee\ga_{\{1\}}$
such that, for $x\in\ga_{\{0\}}$, $y\in\ga_{\{1\}}$,
$$
\ve(xy)=\tr(xp)\f(y)p+\psi(y)p^{\perp}xp^{\perp}\,,
$$
respectively
$$
\ve(xy)=\tr(xp)\f(y)p+\psi(p^{\perp}xp^{\perp}y)p^{\perp}\,.
$$

\subsection{Case 8:} $\cz(\car(\ve))\sim\bc^{2}$, 
$\#$ of orbits of $\Th\lceil_{\cz(\car(\ve))}=1$.

Choose $\{Q_{\pm\chi}\}$ as the generators of $\cz(\car(\ve))$. 
We have two possibilities. First
$$
\ve(x)=\f(Q_{\chi}xQ_{\chi})Q_{\chi}
+\f(Q_{\chi}\Th(x)Q_{\chi})Q_{-\chi}\,,\quad x\in\ga_{[0,1]}\,,
$$
$\f$ being a state on $Q_{\chi}\ga_{[0,1]}Q_{\chi}$. Second, let 
$\gb\subset\ga_{[0,1]}$ be the tensor completion of $\ga_{\{0\}}$ in 
$\ga_{[0,1]}$.\footnote{According to 
\eqref{kw}, the subalgebra $\gb$ is obtained by 
constructing a systems $\{e_{kl}(j)\,,f_{kl}(j)\,|\,j,k,l=1,2\}$  
of four mutually commuting $2\times2$ matrix units for $\ga_{[0,1]}$. 
Notice that $\gb$ is localized in the whole $\ga_{[0,1]}$, 
and is $\Th$--invariant.}
Then there exists a state $\f$ on $\gb$ such that for 
$x\in\ga_{\{0\}}$, $y\in\gb$,
$$
\ve(xy)=\f(y)Q_{\chi}xQ_{\chi}
+\f(\Th(y))Q_{-\chi}xQ_{-\chi})\,.
$$

\subsection{Case 9:} $\cz(\car(\ve))\sim\bc$.

We treat only one possibility, the two remaining ones generating  
one step product states (see e.g. Subsection 6.1 of \cite{AFM}). 
Let $N$, $\bar N$ be the algebra generated by $a_{1}, a_{1}^{+}$,
$a_{2}, a_{2}^{+}$ respectively. Then there exists an even state $\f$ on 
$\bar N\bigvee\ga_{\{1\}}$ such that 
for $x\in N$, $y\in\bar N\bigvee\ga_{\{1\}}$,
$$
\ve(xy)=\f(y)x\,.
$$
Notice that this example is nothing but that the two
block factor 
treated in Subsection 6.3 of \cite{AFM}. This is easily seen by 
passing in \cite{AFM}, to the two point regrouped algebra.

\subsection{Case 10:} two examples with $\ga_{\{n\}}\sim\bm_{2^{3}}(\bc)$.

We describe two examples relative to more complicated situations 
than the previous ones. Let $\{a_{i},a^{+}_{i}\,|\,i=1,2,3\}$,
$\{b_{i},b^{+}_{i}\,|\,i=1,2,3\}$ be the generators of $\ga_{\{0\}}$,
$\ga_{\{1\}}$ respectively. Let $\{e_{kl}(j)\,,f_{kl}(j)\,|\,j,k,l=1,2\}$
of commuting 
$2\times2$ matrix units obtained according to \eqref{kw}, and 
realizing the isomorphism 
$\ga_{[0,1]}\sim\underbrace{\bm_{2}(\bc)\otimes\cdots\otimes\bm_{2}
(\bc)}_{6\text{--times}}$. Put for $\chi\in\bt$,
$$
P_{\chi}:=\frac{1}{2}\big(\idd+\chi 
a_{1}+\overline{\chi}a_{1}^{+}\big)\,.
$$
First define $N_{i}$, $\bar N_{i}$ as the algebras generated by 
$\{e_{ii}(1)a_{2}\,,e_{ii}(1)a^{+}_{2}\}$, 
$\{e_{ii}(1)a_{3}\,,e_{ii}(1)a^{+}_{3}\}$, $i=1,2$ respectively. 
Choose even states $\f_{i}$ on $\bar N_{i}\bigvee\ga_{\{1\}}$. Then 
for $x_{i}\in N_{i}$, $y\in\bar N_{i}\bigvee\ga_{\{1\}}$,
$$
\ve\big(\sum_{i=1}^{2}x_{i}y_{i}\big)=
\sum_{i=1}^{2}\f_{i}(y_{i})x_{i}\,.
$$
Second define $N_{\chi}$, $M_{\chi}$ as the algebras generated by 
$\{P_{\chi}e_{ij}(2)\,|\,i,j=1,2\}$, 
$\{P_{\chi}e_{ij}(3)f_{kl}(n)\,|\,i,j,k,l=1,2\,,n=1,2,3\}$ respectively. 
Choose a state $\f$ on $M_{\chi}$. Then 
for $x_{\pm\chi}\in N_{\pm\chi}$, $y_{\pm\chi}\in M_{\pm\chi}$,
$$
\ve(x_{\chi}y_{\chi}+x_{-\chi}y_{-\chi})=
\f(y_{\chi})x_{\chi}+\f(\Th(y_{-\chi}))x_{-\chi}\,.
$$

\section*{acknowledgements}

The author is grateful to L. Accardi for suggesting the themes and for several useful discussions. He also acknowledges one of the referees whose suggestions considerably contributed to improve the presentation of the present work.

\end{document}